\documentclass[pdflatex,sn-mathphys-num]{sn-jnl}
\usepackage{amsmath,amssymb,bm,mathtools,booktabs,array,microtype}
\usepackage{graphicx,placeins,etoolbox,xcolor}
\hypersetup{hidelinks,
 pdftitle={Minimax strong duality and sharp convex-hull limits for the projected subgradient method},
 pdfsubject={Finite-horizon minimax strong duality, optimal stepsizes, and sharp online averaging constants (Springer Nature, December 2024 template)},
 pdfkeywords={projected subgradient method, convex hull, minimax strong duality, online averaging}}

\makeatletter
\renewcommand\normalsize{%
 \@setfontsize\normalsize{11}{13.5}%
 \abovedisplayskip=9pt plus 2pt minus 2pt%
 \belowdisplayskip=9pt plus 2pt minus 2pt%
 \abovedisplayshortskip=2pt plus 2pt%
 \belowdisplayshortskip=5pt plus 2pt minus 1pt%
 \let\@listi\@listI}
\DeclareMathSizes{11}{11}{8}{6}
\renewcommand\abstractfont{\reset@font\unboldmath\fontsize{10.5}{12.5}\selectfont
 \leftskip=0pt\rightskip=0pt\parfillskip=0pt plus 1fil}
\renewcommand\abstractheadfont{\reset@font\unboldmath\fontsize{10.5}{12.5}\bfseries\selectfont\titraggedcenter}
\renewcommand\keywordfont{\reset@font\unboldmath\fontsize{10}{12}\selectfont
 \leftskip=0pt\rightskip=0pt plus 0.5fil}
\renewcommand\addressfont{\reset@font\fontsize{10.5}{12.5}\selectfont\titraggedcenter}
\renewcommand\figurecaptionfont{\reset@font\fontsize{9.5}{11.5}\selectfont}
\renewcommand\tablecaptionfont{\reset@font\fontsize{9.5}{11.5}\selectfont}
\renewcommand\tablebodyfont{\reset@font\fontsize{10}{12}\selectfont}
\patchcmd{\@maketitle}{\vskip21pt}{\vskip0pt}{}{%
 \PackageWarning{psg-layout}{Title top spacing patch not applied}}
\patchcmd{\@maketitle}{\vskip20pt\nointerlineskip}{\vskip6pt\nointerlineskip}{}{}
\patchcmd{\@maketitle}{\vskip20pt}{\vskip12pt}{}{}
\patchcmd{\@maketitle}{\vskip24pt}{\vskip10pt}{}{}
\patchcmd{\@maketitle}{\vskip24pt}{\vskip12pt}{}{}
\patchcmd{\@maketitle}{\vskip36pt\vskip0pt}{\vskip16pt}{}{}
\def\ps@headings{%
 \let\@oddhead\@empty\let\@evenhead\@empty
 \def\@oddfoot{\hfil\thepage\hfil}\let\@evenfoot\@oddfoot
 \let\@mkboth\markboth}
\let\ps@titlepage\ps@headings
\let\ps@plain\ps@headings
\makeatother
\normalsize
\newtheoremstyle{psgresult}{9pt plus 2pt minus 2pt}{9pt plus 2pt minus 2pt}
 {\itshape}{0pt}{\bfseries}{.}{.5em}{}
\theoremstyle{psgresult}
\newtheorem{theorem}{Theorem}[section]
\newtheorem{lemma}[theorem]{Lemma}
\newtheorem{proposition}[theorem]{Proposition}
\AtBeginEnvironment{proof}{\let\small\normalsize}

\newcommand{\R}{\mathbb R}
\newcommand{\vx}{\bm x}
\newcommand{\vg}{\bm g}
\newcommand{\vu}{\bm u}
\newcommand{\vv}{\bm v}
\newcommand{\vw}{\bm w}
\newcommand{\vz}{\bm z}
\newcommand{\vy}{\bm y}
\newcommand{\vh}{\bm h}
\newcommand{\vn}{\bm n}
\newcommand{\vc}{\bm c}
\newcommand{\vq}{\bm q}
\newcommand{\vb}{\bm b}
\newcommand{\ve}{\bm e}
\newcommand{\vH}{\bm H}
\newcommand{\vX}{\bm X}
\newcommand{\vphi}{\bm\varphi}
\newcommand{\vzeta}{\bm\zeta}
\newcommand{\norm}[1]{\left\|#1\right\|}
\newcommand{\ip}[2]{\left\langle#1,#2\right\rangle}
\newcommand{\cv}{\mathrm{cv}}
\newcommand{\scv}{\mathrm{sc}}
\newcommand{\dd}{\,\mathrm d}
\DeclareMathOperator{\conv}{conv}

\DeclareMathOperator{\proj}{P}
\DeclareMathOperator{\secant}{sec}
\allowdisplaybreaks[2]
\newif\ifauthorinformation
\authorinformationtrue

\begin{document}
\title[Sharp convex-hull limits for the projected subgradient method]{Minimax strong duality and sharp convex-hull limits for the projected subgradient method}
\ifauthorinformation
\author[1]{\fnm{Zhihan} \sur{Zhu}}\email{zhihanzhu@buaa.edu.cn}
\author[2]{\fnm{Yanhao} \sur{Zhang}}\email{yhzhang@cufe.edu.cn}
\author*[1]{\fnm{Yong} \sur{Xia}}\email{yxia@buaa.edu.cn}
\affil[1]{\orgdiv{School of Mathematical Sciences}, \orgname{Beihang University}, \orgaddress{\city{Beijing} \postcode{100191}, \country{China}}}
\affil[2]{\orgdiv{School of National Data Engineering and Security}, \orgname{Central University of Finance and Economics}, \orgaddress{\city{Beijing} \postcode{100098}, \country{China}}}
\hypersetup{pdfauthor={Zhihan Zhu, Yanhao Zhang, Yong Xia}}
\fi
\abstract{The convex hull of a projected subgradient trajectory provides a common
performance benchmark for individual iterates, best observed points,
and weighted averages. For any prescribed stepsizes, we establish
finite-horizon minimax strong duality over admissible instances of
arbitrary dimension: at each horizon one weight vector fixed in advance
achieves the same worst-case error as weights optimized in hindsight for
each instance. We then determine the sharp asymptotic hull limits for two
classical families of diminishing stepsizes. For convex objectives
with feasible points within distance $R$ of a minimizer and queried
subgradients bounded by $G$, the family $aR/(G\sqrt t)$ yields the optimal
leading constant approximately $1.0261$ in units of $RG/\sqrt T$, attained
uniquely at $a\approx0.3303$. For $\mu$-strongly convex objectives,
$a/(\mu t)$ yields the optimal constant $\pi^2/8$ in units of $G^2/(\mu T)$,
attained uniquely at $a=1$. Explicit suffix and sine averages realize
these limits. Cubic and quadratic online averages accommodate an unknown
stopping horizon and have sharp leading constants approximately $0.24\%$
and $0.18\%$ above the respective optimal limits. Experiments on hard
instances exhibit the predicted agreement with the hull optima, and
the proposed averages reduce the number of subgradient queries needed
to reach a target accuracy on ball-constrained sparse regression problems.
The results give precise worst-case benchmarks and practical averaging
rules that nearly attain them with anytime convergence.}
\keywords{Projected subgradient method, Convex hull, Minimax strong duality,
 Stepsize design, Online averaging, Nonsmooth optimization}
\pacs[MSC Classification]{90C25, 90C30, 90C60}
\maketitle

\section{Introduction}\label{sec:introduction}

We consider the nonsmooth convex optimization problem
\begin{equation}\label{eq:problem}
 \min_{\vx\in C} f(\vx),
\end{equation}
where $C\subseteq\R^d$ is nonempty, closed, and convex and
$f:\R^d\to\R$ is convex. We assume a minimizer $\vx^\ast\in C$ exists
and write $f^\ast=f(\vx^\ast)$. All norms are Euclidean. Throughout, the
selected subgradients have norm at most $G>0$. In the convex class, assume
\begin{equation}\label{eq:convex-radius}
 C\subseteq\{\vx:\norm{\vx-\vx^\ast}\le R\},\qquad R>0.
\end{equation}
In the strongly convex class, $f$ is $\mu$-strongly convex for
$\mu>0$ and $C$ may be any nonempty closed convex set.

To solve \eqref{eq:problem}, the projected subgradient
method~\cite{Bubeck2015,Nesterov2018} starts at $\vx_1\in C$
and generates the feasible iterates
\begin{equation}\label{eq:algorithm-intro}
 \vx_{t+1}=\proj_C(\vx_t-\eta_t\vg_t),
 \qquad \vg_t\in\partial f(\vx_t),\qquad \norm{\vg_t}\le G.
\end{equation}
We study two classical families of diminishing stepsizes,
\begin{equation}\label{eq:steps-intro}
 \eta_t=\frac{aR}{G\sqrt t}\quad\text{(convex)},\qquad
 \eta_t=\frac{a}{\mu t}\quad\text{(strongly convex)},
\end{equation}
with a coefficient $a>0$ fixed independently of the stopping horizon.

Convergence guarantees for projected subgradient method typically
specify an output point. Ergodic guarantees concern a prescribed average
of the iterates. Nonergodic guarantees concern individual last iterate. A best-value rule selects an iterate with the smallest observed objective value, with the additional cost of evaluating the objective at each iterate to select the output.

Both families include familiar choices. For convex objectives,
Bubeck~\cite[Sect.~3.1]{Bubeck2015} discusses the schedule
$\eta_t=R/(G\sqrt t)$, corresponding to $a=1$. Under the radius condition,
Zhu, Zhang, and Xia~\cite[Theorem~1]{ZhuZhangXia2025} prove the ergodic bound
\[
 f\left(\frac1T\sum_{t=1}^T\vx_t\right)-f^\ast
 \le \frac{3}{2}\,\frac{RG}{\sqrt T}.
\]
For strongly convex objectives, Beck~\cite[Theorem~8.31]{Beck2017} uses
the schedule in \eqref{eq:steps-intro} with $a=2$. For $T\ge2$, the
classic average with linear weights satisfies
\[
 \bar{\vx}_T^{\mathrm{lin}}
   =\frac{2}{T(T-1)}\sum_{t=1}^T(t-1)\vx_t,\qquad
 f(\bar{\vx}_T^{\mathrm{lin}})-f^\ast\le 2\,\frac{G^2}{\mu T}.
\]
Similar results can also be obtained for the best observed value.

Ergodic, nonergodic, and best-value guarantees each control the error
of the specific point selected by an output rule. All these output points lie in the convex hull of the first $T$ queried
iterates. Define this hull and the simplex of convex-combination weights by
\begin{equation}\label{eq:delta-hull}
 \begin{aligned}
 \mathcal K_T&=\conv\{\vx_1,\ldots,\vx_T\},\\
 \mathcal W_T&=\left\{\vw\in[0,\infty)^T:
                         \sum_{t=1}^T w_t=1\right\}.
 \end{aligned}
\end{equation}
We benchmark the trajectory by the smallest gap in this hull,
\begin{equation}\label{eq:hull-error-intro}
 \Delta_T=\min_{\vw\in\mathcal W_T}
 \left\{ f\left(\sum_{t=1}^T w_t\vx_t\right)-f^\ast \right\}
 =\min_{\vx\in\mathcal K_T} f(\vx)-f^\ast .
\end{equation}
For each stepsize family, we take the worst-case value of
\eqref{eq:hull-error-intro} over the admissible instances for a given $a$
and minimize the resulting asymptotic leading constant over $a>0$.
Let $\mathcal I_T(a)$ denote the admissible instances in the relevant
class at coefficient $a$, including the dimension, function, feasible
set, initial point, and oracle selections. For $\vw\in\mathcal W_T$, let
$E_T(\mathcal I,a,\vw)$ be the objective gap of
$\sum_{t=1}^T w_t\vx_t$ divided by $RG/\sqrt T$ in the convex class or
$G^2/(\mu T)$ in the strongly convex class.

\begin{enumerate}
\item \textbf{Finite-horizon minimax strong duality.}
For every $T$ and every fixed coefficient $a>0$, in either class,
Theorem~\ref{thm:minimax-duality} establishes
\begin{equation}\label{eq:weight-minimax-duality}
 \sup_{\mathcal I\in\mathcal I_T(a)}
       \min_{\vw\in\mathcal W_T} E_T(\mathcal I,a,\vw)
 =\min_{\vw\in\mathcal W_T}
       \sup_{\mathcal I\in\mathcal I_T(a)} E_T(\mathcal I,a,\vw).
\end{equation}
The $\sup\min$ order allows
the weights to depend on the objective and trajectory of each instance.
The $\min\sup$ order uses one weight vector for the whole class at each
horizon $T$. The equality follows from the closure of the admissible
class under orthogonal products with scaled coordinates. It holds for
any prescribed stepsize sequence common to all instances, without
knowing the optimal stepsize or the sharp constants.

\item \textbf{Sharp hull limits and explicit attaining averages.}
For convex objectives, the unique optimal coefficient is
$a_\cv\approx0.3303$, giving the sharp leading constant
$C_\cv\approx1.0261$ in units of $RG/\sqrt T$
(Theorem~\ref{thm:convex}).
For strongly convex objectives, $a=1$ is uniquely optimal and yields
$C_\scv=\pi^2/8$ in units of $G^2/(\mu T)$
(Theorem~\ref{thm:strong}).
Explicit suffix and sine averaging weights, depending only on $t$ and
$T$, attain these constants. Matching lower bounds hold at every
point of the respective hard-instance hulls. Consequently, both
sides of \eqref{eq:weight-minimax-duality} converge to the sharp
constant of the corresponding class at its optimal coefficient.

\item \textbf{Sharp guarantees for near-optimal online averages.}
The weights attaining the sharp constants are set for a prescribed horizon
$T$. To allow an unknown stopping horizon, we use weights $t^3$ with
$a=1/3$ in the convex case and $t^2$ with $a=1$ in the strongly convex case.
The weights depend only on $t$, so the average can be updated as each
iterate arrives without evaluating the objective. This gives anytime
convergence with exact worst-case leading constants $C_3\approx1.0285$
and $C_2=9\log3/8\approx1.2359$, respectively. Their relative excesses
over the optimized hull limits are approximately $0.24\%$ and $0.18\%$
(Theorem~\ref{thm:online}).
\end{enumerate}

\subsection{Related work and perspective}\label{sec:related}

Zhu, Zhang, and Xia~\cite{ZhuZhangXia2025} study weak ergodic averages with
weights $\eta_t^{-k}$. The exponents $k=6$ for square-root steps and $k=2$
for harmonic steps give the cubic and quadratic weights used here.
Our analysis compares these outputs with the sharp hull benchmark.
Nesterov's predefined primal steps~\cite{Nesterov2024} and the
Lipschitz-free schedules of Xia, Zhang, and Zhu~\cite{XiaZhangZhu2025}
address complementary aspects of stepsize selection and constraint-aware
scaling.

For strongly convex stochastic optimization, suffix averaging
\cite{RakhlinShamirSridharan2012} and polynomially weighted averaging
\cite{ShamirZhang2013} recover the optimal order. Lacoste-Julien, Schmidt,
and Bach~\cite{LacosteJulienSchmidtBach2012} analyze linear weights with a
shifted harmonic stepsize.
Harvey, Liaw, Plan, and Randhawa~\cite{HarveyEtAl2019} identify sharp
logarithmic effects in last-iterate behavior. Grimmer and Li
\cite{GrimmerLi2025} develop primal-dual descriptions and convergence
certificates for strongly convex subgradient methods. Our strongly convex
upper bound uses a rank-two certificate that retains the quadratic
contribution of the average.


Exact worst-case analysis and constructive lower functions have a central
role in performance estimation
\cite{DroriTeboulle2014,TaylorHendrickxGlineur2017} and optimal method design
\cite{DroriTaylor2020}. Taylor, Hendrickx, and Glineur
\cite{TaylorHendrickxGlineurComposite2017} extend performance estimation to
composite convex problems, including exact guarantees for the last iterate.
Luner and Grimmer~\cite{LunerGrimmer2026}
show that, in smooth convex optimization, averaging cannot improve the worst-case performance of the last iterate of gradient descent for certain stepsizes.
Zamani and Glineur~\cite{ZamaniGlineur2025} derive
exact last-iterate results for subgradient methods using moving reference points and design an optimal horizon-dependent schedule. Zoll and Grimmer~\cite{ZollGrimmer2026} characterize minimax-optimal fixed-step methods for Lipschitz convex minimization under an initial-distance bound. For the two prescribed update shapes in \eqref{eq:steps-intro}, we use the minimum over the entire iterate hull as the benchmark and connect this worst-case value to stepsize design and instance-independent averaging through minimax strong duality.

\section{The convex-hull benchmark and the main results}\label{sec:main}

Under the assumptions and stepsize families of Section~\ref{sec:introduction},
let $\mathcal I_T^\cv(a)$ and $\mathcal I_T^\scv(a)$ denote the convex
and strongly convex instance classes, respectively. Their parameters
$R,G$ (or $\mu,G$) are common to the instances in each class.
For the hull gap $\Delta_T$ in \eqref{eq:hull-error-intro}, define
\begin{align}
 H_T^\cv(a)&=\sup_{\mathcal I_T^\cv(a)}\frac{\Delta_T}{RG},
 &\Gamma_\cv(a)&=\limsup_{T\to\infty}\sqrt T\,H_T^\cv(a),
 \label{eq:gamma-convex}\\
 H_T^\scv(a)&=\sup_{\mathcal I_T^\scv(a)}\frac{\mu\Delta_T}{G^2},
 &\Gamma_\scv(a)&=\limsup_{T\to\infty}T H_T^\scv(a).
 \label{eq:gamma-strong}
\end{align}
The dimension is unrestricted, and the instance may vary with $T$.
The coefficient $a$ is fixed before the asymptotic limit is taken.
All upper remainders below are uniform over the relevant instance class.

\begin{theorem}[Finite-horizon minimax strong duality]\label{thm:minimax-duality}
For every $T\ge1$ and $a>0$, in either instance class,
\begin{equation}\label{eq:minimax-limit}
 \underbrace{\sup_{\mathcal I\in\mathcal I_T(a)}
       \min_{\vw\in\mathcal W_T} E_T(\mathcal I,a,\vw)}_{L_T(a)}
 =\underbrace{\min_{\vw\in\mathcal W_T}
       \sup_{\mathcal I\in\mathcal I_T(a)}
       E_T(\mathcal I,a,\vw)}_{U_T(a)}.
\end{equation}
In particular, a weight vector chosen in advance for the entire
instance class attains the finite-horizon worst-case hull gap.
The identity also holds for any common prescribed positive stepsize
sequence in either class.
\end{theorem}
\begin{proof}
Fix $T$, $a$, and the corresponding prescribed stepsizes, which are
common to all instances. For $\vw\in\mathcal W_T$, write
$e_{\mathcal I}(\vw)=E_T(\mathcal I,a,\vw)$. Each $e_{\mathcal I}$ is
continuous and convex in $\vw$, since the objective is finite and convex
and the weighted output is affine in $\vw$. In the convex class the
subgradient inequality gives $f(\vx_t)-f^\ast\le RG$. In the strongly
convex class it gives
\[
 f(\vx_t)-f^\ast\le
 G\norm{\vx_t-\vx^\ast}-\frac{\mu}{2}\norm{\vx_t-\vx^\ast}^2
 \le\frac{G^2}{2\mu}.
\]
By convexity, the same bounds hold for every convex combination of the
iterates. Hence $\sup_{\mathcal I\in\mathcal I_T(a)}e_{\mathcal I}(\vw)$
is finite and lower semicontinuous. Compactness of $\mathcal W_T$
ensures attainment of both the minimum of each $e_{\mathcal I}$ and the
minimum defining $U_T(a)$. The pointwise inequality
\[
 \min_{\widetilde{\vw}\in\mathcal W_T}e_{\mathcal I}(\widetilde{\vw})
 \le e_{\mathcal I}(\vw)
\]
implies the weak duality inequality $L_T(a)\le U_T(a)$ after taking
the supremum over instances and then the minimum over $\vw$.

To prove the reverse inequality, it suffices to show that for every
$c<U_T(a)$ there is a single admissible instance whose normalized hull
gap exceeds $c$. The definition of $U_T(a)$ provides a difficult instance
for each individual weight vector. To obtain one instance that is
difficult for all weight vectors simultaneously, we first establish
closure under finite convex combinations of the error functions.

\emph{Closure under finite convex combinations.}
Consider $k$ admissible instances $\mathcal I_i\in\mathcal I_T(a)$,
$i=1,\ldots,k$, with dimensions $d_i$, objectives
$f_i:\R^{d_i}\to\R$, feasible sets $C_i\subseteq\R^{d_i}$, and
initial points $\vx_{i,1}\in C_i$. Let $\vx_i^\ast$ be a minimizer
satisfying the radius condition in the convex class, and set
$f_i^\ast=f_i(\vx_i^\ast)$. Denote by $\vx_{i,t}$ the $t$th
iterate of instance $\mathcal I_i$ and by $\vg_{i,t}$ its selected
subgradient. All these trajectories use the same stepsizes:
\[
 \vx_{i,t+1}=\proj_{C_i}(\vx_{i,t}-\eta_t\vg_{i,t}),
 \qquad \vg_{i,t}\in\partial f_i(\vx_{i,t}).
\]
Choose coefficients $p_i\ge0$, $\sum_{i=1}^k p_i=1$, for combining
the instances. Unlike $w_t$, which averages iterates within an instance,
$p_i$ weights the $i$th instance. We seek one admissible instance
$\mathcal I_p$ realizing the error function
$\sum_{i=1}^k p_i e_{\mathcal I_i}$.
Omit zero-coefficient blocks and reindex the remaining blocks so that
$p_i>0$ for all $i$.

In the product space $\R^{d_1}\times\cdots\times\R^{d_k}\cong\R^D$,
where $D=\sum_{i=1}^k d_i$, write $\vy=(\vy^{(1)},\ldots,\vy^{(k)})$
with $\vy^{(i)}\in\R^{d_i}$ denoting the $i$th block. Define
\[
 F(\vy)=\sum_{i=1}^k p_i f_i(\vy^{(i)}/\sqrt{p_i}),
 \qquad C_p=\prod_{i=1}^k\sqrt{p_i}\,C_i .
\]
Here $\sqrt{p_i}\,C_i=\{\sqrt{p_i}\vx:\vx\in C_i\}$, and the
product is the Cartesian product of these scaled feasible sets.
The combined instance starts from the block vector
\[
 \vy_1=(\sqrt{p_1}\vx_{1,1},\ldots,\sqrt{p_k}\vx_{k,1}).
\]
The objective $F$ is finite and convex, and $C_p$ is nonempty, closed,
and convex. Its minimum is $F^\ast=\sum_{i=1}^k p_i f_i^\ast$,
attained at $\vy^\ast=(\sqrt{p_i}\vx_i^\ast)_{i=1}^k$.

We now show that the trajectory of the combined instance consists of
the scaled original trajectories in their respective blocks.
Suppose $\vy_t=(\sqrt{p_i}\vx_{i,t})_{i=1}^k$. By the subgradient
scaling rule and separability of $F$ we can select
$\vh_t=(\sqrt{p_i}\vg_{i,t})_{i=1}^k\in\partial F(\vy_t)$.
Euclidean projection onto $C_p$ separates into the blocks, and
$\proj_{\sqrt{p_i}C_i}(\sqrt{p_i}\vv)=\sqrt{p_i}\proj_{C_i}(\vv)$
for $\vv\in\R^{d_i}$. Therefore
\[
 \begin{aligned}
 \vy_{t+1}&=\proj_{C_p}(\vy_t-\eta_t\vh_t)\\
 &=\bigl(\sqrt{p_i}\proj_{C_i}(\vx_{i,t}-\eta_t\vg_{i,t})\bigr)_{i=1}^k\\
 &= (\sqrt{p_i}\vx_{i,t+1})_{i=1}^k.
 \end{aligned}
\]
Starting from $\vy_1$, induction gives
$\vy_t=(\sqrt{p_i}\vx_{i,t})_{i=1}^k$ for every $1\le t\le T$.
Moreover, $\norm{\vh_t}^2=\sum_{i=1}^k p_i\norm{\vg_{i,t}}^2\le G^2$.
In the convex class, any feasible point
$\vy=(\sqrt{p_i}\vu_i)_{i=1}^k\in C_p$, with $\vu_i\in C_i$, satisfies
\[
 \norm{\vy-\vy^\ast}^2=\sum_{i=1}^k p_i\norm{\vu_i-\vx_i^\ast}^2
 \le R^2.
\]
In the strongly convex class, writing
$\vu_i=\vy^{(i)}/\sqrt{p_i}$, each scaled block function satisfies
\[
 p_i f_i\left(\frac{\vy^{(i)}}{\sqrt{p_i}}\right)
 -\frac{\mu}{2}\norm{\vy^{(i)}}^2
 =p_i\left(f_i(\vu_i)-\frac{\mu}{2}\norm{\vu_i}^2\right).
\]
The right-hand side is convex in $\vy^{(i)}$. Summing these functions over the
orthogonal blocks shows that $F$ is $\mu$-strongly convex. This
explains the use of $\sqrt{p_i}$ in the coordinate scaling: it
preserves the quadratic norm bounds and strong convexity parameter.
Thus $\mathcal I_p$ belongs to the same admissible instance class.

For $\vw\in\mathcal W_T$, write
$\bar{\vx}_i(\vw)=\sum_{t=1}^T w_t\vx_{i,t}$ and
$\bar{\vy}(\vw)=\sum_{t=1}^T w_t\vy_t$ for the outputs of the
original and combined instances, respectively. The block representation
gives $\bar{\vy}(\vw)=(\sqrt{p_i}\bar{\vx}_i(\vw))_{i=1}^k$, so
\[
 F(\bar{\vy}(\vw))-F^\ast
 =\sum_{i=1}^k p_i\bigl(f_i(\bar{\vx}_i(\vw))-f_i^\ast\bigr).
\]
Dividing by the common normalizing scale gives the exact equality
\begin{equation}\label{eq:product-loss}
 E_T(\mathcal I_p,a,\vw)=\sum_{i=1}^k p_i E_T(\mathcal I_i,a,\vw)
 \qquad(\vw\in\mathcal W_T).
\end{equation}
Thus any finite convex combination of admissible error functions is
realized by one admissible instance.

\emph{The reverse minimax inequality.}
Fix any $c<U_T(a)$. For every $\vw\in\mathcal W_T$, some
$\mathcal I\in\mathcal I_T(a)$ satisfies $e_{\mathcal I}(\vw)>c$.
Since $e_{\mathcal I}$ is continuous, the relatively open sets
$\{\vw\in\mathcal W_T:e_{\mathcal I}(\vw)>c\}$ cover the compact
simplex. Choose a finite subcover corresponding to the instances
$\mathcal I_1,\ldots,\mathcal I_k$, and write
$e_i(\vw)=e_{\mathcal I_i}(\vw)$. We have
\[
 \max_{1\le i\le k}e_i(\vw)>c\qquad(\vw\in\mathcal W_T).
\]
The remaining task is to find a single vector of combination
coefficients $\bm p$ for these instances such that
$\sum_{i=1}^k p_i e_i(\vw)>c$ for all $\vw$. To this end, define the
upper image
\[
 \mathcal A=\{\vz\in\R^k:\exists\,\vw\in\mathcal W_T
                 \text{ with } e_i(\vw)\le z_i\ (1\le i\le k)\}.
\]
A vector $\vz$ belongs to $\mathcal A$ if one common averaging vector
achieves errors no larger than the corresponding coordinates of $\vz$.
This set is closed: for a convergent sequence in $\mathcal A$, a
convergent subsequence of the associated weights exists by compactness
of $\mathcal W_T$, and continuity of the $e_i$ gives a witnessing
weight at the limit. It is convex because convex combinations of
witnessing weights satisfy the corresponding upper bounds by convexity
of each $e_i$.

The preceding strict inequality implies $c\mathbf 1_k\notin\mathcal A$,
where $\mathbf 1_k\in\R^k$ is the vector of ones.
Strict separation yields a nonzero $\vq\in\R^k$ satisfying
\[
 c\sum_{i=1}^k q_i < \inf_{\vz\in\mathcal A}\ip{\vq}{\vz}.
\]
Since $\mathcal A+[0,\infty)^k=\mathcal A$, a negative component of
$\vq$ would make this infimum $-\infty$. Hence $\vq\ge0$ and we can
normalize it to $\bm p=\vq/\sum_{i=1}^k q_i$. Each $\vz\in\mathcal A$
dominates the error vector of some $\vw$, and every such error vector
belongs to $\mathcal A$. Since $\bm p\ge0$,
\[
 c<\inf_{\vz\in\mathcal A}\ip{\bm p}{\vz}
 =\min_{\vw\in\mathcal W_T}\sum_{i=1}^k p_i e_i(\vw).
\]
Applying \eqref{eq:product-loss} with these fixed coefficients yields
one admissible instance $\mathcal I_p$ satisfying
\[
 c<\min_{\vw\in\mathcal W_T} E_T(\mathcal I_p,a,\vw)
   \le L_T(a).
\]
Letting $c\uparrow U_T(a)$ proves $U_T(a)\le L_T(a)$.
Together with weak duality this gives \eqref{eq:minimax-limit}.
Because the minimum defining $U_T(a)$ is attained, a single weight
vector fixed in advance achieves the common worst-case value over the
entire class. The product construction uses unrestricted dimension
and common prescribed stepsizes, but no particular stepsize formula.
The identity therefore holds for any common prescribed positive
stepsize sequence in either class. 
\end{proof}

We next determine the best asymptotic values of this common benchmark
and give explicit weights that attain them.

\begin{theorem}[Sharp convex hull limit]\label{thm:convex}
Let $\tau_\ast=0.575692574965\ldots$ be the unique root in $(1/3,1)$ of
\begin{equation}\label{eq:tau-root}
 \log\frac{2(1+\tau+2\tau^2)}{(1-\tau)^2}
 =\frac{1+3\tau}{2(1-\tau)}.
\end{equation}
Set
\begin{equation}\label{eq:convex-constants}
 \begin{aligned}
 C_\cv&=\frac1{\sqrt{(1-\tau_\ast)(1+\tau_\ast+2\tau_\ast^2)}}
           =1.026071437300\ldots,\\
 a_\cv&=\sqrt{\frac{\tau_\ast(1-\tau_\ast)}
                         {1+\tau_\ast+2\tau_\ast^2}}
           =0.330334203702\ldots.
 \end{aligned}
\end{equation}
Then $\inf_{a>0}\Gamma_\cv(a)=C_\cv$, with a unique minimizing coefficient
$a=a_\cv$. At this coefficient,
$\lim_{T\to\infty}\sqrt T H_T^\cv(a_\cv)=C_\cv$.
For $n=\lceil\tau_\ast T\rceil$ and sufficiently large $T$, define
\begin{equation}\label{eq:convex-optimal-weights}
 w_{t,T}=\begin{cases}
 0,&t<n,\\[2pt]
 \displaystyle\frac{(T-t+1)\sqrt{t-1}-(T-t)\sqrt t}
 {(T-n+1)\sqrt{n-1}},&n\le t\le T.
 \end{cases}
\end{equation}
These weights are nonnegative, sum to one, and increase on their support.
With $a=a_\cv$,
\begin{equation}\label{eq:convex-attainment}
 f\left(\sum_{t=1}^T w_{t,T}\vx_t\right)-f^\ast
 \le\sum_{t=1}^T w_{t,T}\bigl(f(\vx_t)-f^\ast\bigr)
 \le (C_\cv+o(1))\frac{RG}{\sqrt T}.
\end{equation}
\end{theorem}

\begin{theorem}[Sharp strongly convex hull limit]\label{thm:strong}
For every fixed $a>0$,
\begin{equation}\label{eq:strong-universal-lower}
 \liminf_{T\to\infty}T H_T^\scv(a)\ge\frac{\pi^2}{8}.
\end{equation}
For $a=1$, the average
\begin{equation}\label{eq:sine-average}
 \vx_T^{\sin}=
 \frac{\sum_{t=1}^T t\sin(\pi t/(2T))\vx_t}
      {\sum_{t=1}^T t\sin(\pi t/(2T))}
\end{equation}
satisfies
\begin{equation}\label{eq:sine-bound}
 f(\vx_T^{\sin})-f^\ast
 \le\frac{G^2}{\mu}
      \left(\frac{\pi^2}{8T}+O(T^{-3/2})\right).
\end{equation}
The inequality in \eqref{eq:strong-universal-lower} is strict
for every $a\ne1$. Consequently,
\begin{equation}\label{eq:strong-optimum}
 \inf_{a>0}\Gamma_\scv(a)=\Gamma_\scv(1)
 =C_\scv:=\frac{\pi^2}{8}=1.233700550136\ldots.
\end{equation}
The coefficient $a=1$ is uniquely optimal, and
$\lim_{T\to\infty}T H_T^\scv(1)=C_\scv$.
\end{theorem}

Theorems~\ref{thm:convex} and \ref{thm:strong} are proved in
Sections~\ref{sec:convex} and \ref{sec:strong}, respectively.
Together with Theorem~\ref{thm:minimax-duality}, they imply the sharp
asymptotic realization
\begin{equation}\label{eq:asymptotic-duality}
 \lim_{T\to\infty}L_T(a_\star)
 =\lim_{T\to\infty}U_T(a_\star)=\lim_{T\to\infty}\sup_{\mathcal I\in\mathcal I_T(a_\star)} E_T(\mathcal I,a_\star,\vw_T^\star)=C_\star.
\end{equation}
Here $(a_\star,C_\star)=(a_\cv,C_\cv)$ or $(1,\pi^2/8)$, and
$\vw_T^\star$ consists of the weights in
\eqref{eq:convex-optimal-weights} or \eqref{eq:sine-average}.
Indeed, $L_T(a_\star)$ lower bounds the worst-case error of these
weights, while \eqref{eq:convex-attainment} and \eqref{eq:sine-bound}
give matching upper bounds.

\begin{theorem}[Sharp constants for near-optimal online averages]\label{thm:online}
In the convex class, use $\eta_t=R/(3G\sqrt t)$ and return
\begin{equation}\label{eq:cubic-average}
 \vx_T^{(3)}=\frac{\sum_{t=1}^T t^3\vx_t}{\sum_{t=1}^T t^3}.
\end{equation}
Let $u=0.161616224683\ldots$ be the root in
$(0.16161622,0.16161623)$ of $16-99u+83u^9=0$, and put
\[
 L(s)=\sum_{j=0}^7s^j,\qquad P(s)=\sum_{j=0}^7(8-j)s^j.
\]
Define
\begin{equation}\label{eq:C3-definition}
 C_3=\frac1{12u}+\frac{3u^7}{28}
            +\frac34\int_u^1\frac{L(s)^2}{sP(s)}\dd s=1.028520914880\ldots .
\end{equation}
Then
\begin{equation}\label{eq:cubic-main}
 \lim_{T\to\infty}\sqrt T
 \sup_{\mathcal I_T^\cv(1/3)}
 \frac{f(\vx_T^{(3)})-f^\ast}{RG}= C_3.
\end{equation}
In the strongly convex class, use $\eta_t=1/(\mu t)$ and return
\begin{equation}\label{eq:quadratic-average}
 \vx_T^{(2)}=\frac{\sum_{t=1}^T t^2\vx_t}{\sum_{t=1}^T t^2}.
\end{equation}
The uniform upper bound
\begin{equation}\label{eq:quadratic-main}
 f(\vx_T^{(2)})-f^\ast
 \le\frac{G^2}{\mu}
 \left(\frac{9\log3}{8T}+O(T^{-3/2})\right)
\end{equation}
has the sharp leading constant
\begin{equation}\label{eq:quadratic-sharp}
 \lim_{T\to\infty}T\sup_{\mathcal I_T^\scv(1)}
     \frac{\mu(f(\vx_T^{(2)})-f^\ast)}{G^2}
 =C_2:=\frac{9\log3}{8}=1.235938824751\ldots .
\end{equation}
The averages admit the online recursions, initialized at $\vx_1$,
\begin{align}
 \vx_t^{(3)}&=\vx_{t-1}^{(3)}
       +\frac{4t}{(t+1)^2}(\vx_t-\vx_{t-1}^{(3)}),
       \label{eq:cubic-recursion}\\
 \vx_t^{(2)}&=\vx_{t-1}^{(2)}
       +\frac{6t}{(t+1)(2t+1)}(\vx_t-\vx_{t-1}^{(2)}).
       \label{eq:quadratic-recursion}
\end{align}
\end{theorem}

Theorem~\ref{thm:online} is proved in Section~\ref{sec:online}.
Table~\ref{tab:main} summarizes the benchmark and its online realization. 

\begin{table}[!htbp]
\caption{Sharp hull benchmarks and online realizations. In each column,
the error coefficient multiplies the scale in the first row.}
\label{tab:main}
\centering
\renewcommand{\arraystretch}{1.12}
\begin{tabular}{@{}>{\raggedright\arraybackslash}p{.36\textwidth}>{\raggedright\arraybackslash}p{.29\textwidth}>{\raggedright\arraybackslash}p{.29\textwidth}@{}}
\toprule
 & Convex & Strongly convex\\
\midrule
Error scale & $RG/\sqrt T$ & $G^2/(\mu T)$\\
Sharp hull coefficient & $C_\cv\approx1.0261$ & $\pi^2/8\approx1.2337$\\
Unique optimal step coefficient & $a_\cv\approx0.3303$ & $a=1$\\
Attaining output & Increasing suffix, \eqref{eq:convex-optimal-weights}
 & Sine weights, \eqref{eq:sine-average}\\
Online step coefficient & $a=1/3$ & $a=1$\\
Online output weights & $t^3$ & $t^2$\\
Sharp online coefficient & $C_3\approx1.0285$ & $9\log3/8\approx1.2359$\\
Relative excess & $\approx0.24\%$ & $\approx0.18\%$\\
\bottomrule
\end{tabular}
\end{table}

The exact relative excesses, expressed in percent, are
\[
 100(C_3/C_\cv-1)\%\approx0.24\%,\qquad
 100(9\log3/\pi^2-1)\%\approx0.18\%.
\]
These compare worst-case leading coefficients. The recursions themselves
use neither the horizon nor function values. They average the entire history
and add one vector to the storage needed for the current iterate.

\FloatBarrier

\section{The sharp convex hull benchmark}\label{sec:convex}

Under the change of variables $\vy=(\vx-\vx^\ast)/R$, replace $C$ by
$(C-\vx^\ast)/R$ and $f$ by $[f(\vx^\ast+R\vy)-f^\ast]/(RG)$.
Euclidean projection commutes with this transformation, the selected
subgradients are divided by $G$, and the stepsizes become $a/\sqrt t$.
We may therefore assume throughout this section that $R=G=1$,
$\vx^\ast=0$, and $f^\ast=0$, so that $0\in C$, $\norm{\vx}\le1$
on $C$, and $\eta_t=a/\sqrt t$.
The proof couples a weighted projection potential to a family of trajectories
with a common active supporting plane.

\subsection{A projection potential and an attaining average}

\begin{lemma}[Weighted projection potential]\label{lem:cv-potential}
Fix a segment of iteration indices $1\le m\le t\le n\le T$, with
nonnegative averaging weights $\omega_t$ such that
$\sum_{t=m}^n\omega_t=1$. Choose auxiliary scalars
$1=S_{m-1}\le S_m\le\cdots\le S_n$.
Write $W_t=\sum_{j=t}^n\omega_j$ for the weight remaining from time $t$
onward, set $W_{n+1}=0$, and define
\begin{equation}\label{eq:prescribed-coefficients}
 \begin{aligned}
 \lambda_t&=\eta_t S_t,\quad \rho_t=S_t-S_{t-1},\\
 \nu_t&=\frac{W_t}{S_{t-1}}-\frac{W_{t+1}}{S_t},
 \qquad \kappa_t=\frac{\nu_t}{\lambda_t}.
 \end{aligned}
\end{equation}
Then, with $[x]_+=\max\{x,0\}$,
\begin{align}
 \sum_{t=m}^n\omega_t f(\vx_t)
 &\le\frac12\left\{\kappa_m+
 \sum_{t=m+1}^n[\kappa_t-\kappa_{t-1}]_+(2S_{t-1}-1)^2
 \right\}+\frac12\sum_{t=m}^n\nu_t\lambda_t .
 \label{eq:prescribed-bound}
\end{align}
In particular, write $M_t=\sum_{j=t}^n\lambda_j$, $M_{n+1}=0$,
and $\kappa=M_m^{-1}$. If the induced weights
\begin{equation}\label{eq:cv-induced-weights}
 w_t=\kappa\{S_{t-1}\lambda_t-(S_t-S_{t-1})M_{t+1}\}
\end{equation}
are nonnegative, they sum to one and satisfy
\begin{equation}\label{eq:cv-potential-bound}
 \sum_{t=m}^n w_t f(\vx_t)
 \le\frac{1+\sum_{t=m}^n\lambda_t^2}{2M_m}.
\end{equation}
\end{lemma}
\begin{proof}
All sums in this proof are over the segment unless specified otherwise.
The scalars $\nu_t$ and $\kappa_t$ are nonnegative since the tails $W_t$
are nonincreasing and $S_t$ is nondecreasing. Telescoping their tails gives
$\sum_{j>t}\nu_j=W_{t+1}/S_t$, so
\[
 \omega_t=S_{t-1}\nu_t-\rho_t\sum_{j>t}\nu_j.
\]
To apply the projection inequality with a moving reference point,
introduce
\[
 \vu_t=S_{t-1}\vx_t-\sum_{i=m}^{t-1}\rho_i\vx_i,
 \qquad
 \vz_t=S_t^{-1}\sum_{i=m}^t\rho_i\vx_i\in C.
\]
The coefficients of the iterates in $\vz_t$ sum to $(S_t-1)/S_t$,
with the remaining weight on $0\in C$. Thus $\vz_t\in C$.
Since $\vu_t=S_t(\vx_t-\vz_t)$ and
$\vu_{t+1}=S_t(\vx_{t+1}-\vz_t)$, the projection inequality gives
\begin{align*}
 \norm{\vu_{t+1}}^2
 &\le S_t^2\norm{\vx_t-\eta_t\vg_t-\vz_t}^2\\
 &=\norm{\vu_t}^2-2\lambda_t\ip{\vg_t}{\vu_t}
                          +\lambda_t^2\norm{\vg_t}^2 .
\end{align*}
Using $\norm{\vg_t}\le1$ and rearranging yields
\begin{equation}\label{eq:cv-projection-identity}
 \lambda_t\ip{\vg_t}{\vu_t}
 \le\frac{\norm{\vu_t}^2-\norm{\vu_{t+1}}^2}{2}
       +\frac{\lambda_t^2}{2}.
\end{equation}
We next relate the potential decrease to the weighted objective values.
From $\vu_t=\vx_t+\sum_{i=m}^{t-1}\rho_i(\vx_t-\vx_i)$ and the
subgradient inequalities at $0$ and $\vx_i$, we obtain
\[
 \ip{\vg_t}{\vu_t}\ge S_{t-1}f(\vx_t)
                     -\sum_{i=m}^{t-1}\rho_i f(\vx_i).
\]
Multiplying by $\nu_t$ and interchanging the finite sums yields
\[
 \begin{aligned}
 \sum_{t=m}^n\nu_t\ip{\vg_t}{\vu_t}
 &\ge\sum_{t=m}^n
 \left(S_{t-1}\nu_t-\rho_t\sum_{j>t}\nu_j\right)f(\vx_t)\\
 &=\sum_{t=m}^n\omega_t f(\vx_t).
 \end{aligned}
\]
Multiply \eqref{eq:cv-projection-identity} by $\kappa_t$ and sum.
The resulting potential difference is
\[
 \kappa_m\norm{\vu_m}^2+
 \sum_{t=m+1}^n(\kappa_t-\kappa_{t-1})\norm{\vu_t}^2
 -\kappa_n\norm{\vu_{n+1}}^2.
\]
The triangle inequality gives $\norm{\vu_m}\le1$ and
$\norm{\vu_t}\le S_{t-1}+(S_{t-1}-1)=2S_{t-1}-1$.
Dropping the negative terms proves \eqref{eq:prescribed-bound}.

For \eqref{eq:cv-induced-weights}, note that
$w_t=\kappa(S_{t-1}M_t-S_tM_{t+1})$. Their tails are
$W_t=\kappa S_{t-1}M_t$, with $W_m=1$. Substitution into
\eqref{eq:prescribed-coefficients} gives $\nu_t=\kappa\lambda_t$ and
$\kappa_t=\kappa$ throughout the segment. Equation
\eqref{eq:prescribed-bound} reduces to \eqref{eq:cv-potential-bound}.
\end{proof}

\begin{proposition}\label{prop:cv-upper}
Fix two cutoff parameters $0<u<v<1$ with $v^2>1/3$. Write
\begin{equation}\label{eq:cv-dhI}
 d=3v+v^{-1},\quad h=d-4u,\quad
 I(u,v)=\frac1v+\frac2d\log\frac{v(d-4u)}{u(d-4v)}.
\end{equation}
For the weights \eqref{eq:convex-optimal-weights} with
$n=\lceil v^2T\rceil$, every fixed $a>0$ satisfies
\begin{equation}\label{eq:cv-upper-functional}
 \limsup_{T\to\infty}\sqrt T\sup_{\mathcal I_T^\cv(a)}
       \sum_t w_{t,T}f(\vx_t)
 \le U(a,u,v):=\frac1{2ah}+\frac a2I(u,v).
\end{equation}
\end{proposition}
\begin{proof}
Let $m=\lceil u^2T\rceil$, $n=\lceil v^2T\rceil$, and $N=T-n+1$
be the number of suffix iterates. For large $T$, $2\le m<n$.
We apply Lemma~\ref{lem:cv-potential} to $m,\ldots,T$ and design
its auxiliary scalars so that the induced weights vanish before $n$
and $\lambda_t$ is constant from $n$ onward. Define a sequence $p_t$
backwards by
\[
 \kappa_0=\frac{\eta_{n-1}}N,\qquad p_n=N,\qquad
 p_t=p_{t+1}+\frac{\eta_t}{\kappa_0p_{t+1}}
 \quad(t=n-1,\ldots,m).
\]
Set $\ell=\kappa_0p_m$, $\kappa=\kappa_0/\ell^2$, and
\[
 M_t=\begin{cases}
 \ell p_t,&m\le t<n,\\
 \ell(T-t+1),&n\le t\le T+1.
 \end{cases}
\]
Take $\lambda_t=M_t-M_{t+1}$, $S_t=\lambda_t/\eta_t$, and $S_{m-1}=1$.
Then $\kappa M_m=1$. On $m\le t<n$, the backward recursion gives
$\lambda_t=\ell\eta_t/(\kappa_0p_{t+1})$ and hence
$S_t=(\kappa M_{t+1})^{-1}$. Since $S_{m-1}=1=(\kappa M_m)^{-1}$,
the sequence increases on this segment. The weight tails
$\kappa S_{t-1}M_t$ are one for $m\le t\le n$, so the induced
weights vanish on $m,\ldots,n-1$. Extend them by zero to $t<m$.
On the suffix, $S_t=\ell\sqrt t/a$, which joins the middle segment
monotonically since $S_{n-1}=\ell/\eta_{n-1}$. Using
$\kappa\ell^2/a=1/(N\sqrt{n-1})$ in
\eqref{eq:cv-induced-weights} gives \eqref{eq:convex-optimal-weights}.
Its numerator equals
\[
 \sqrt{t-1}-\frac{T-t}{\sqrt t+\sqrt{t-1}},
\]
which increases in $t$ and is positive for $t\ge n$ when $T$ is large.

It remains to evaluate the two scalar terms of
\eqref{eq:cv-potential-bound}. For the limiting calculation, put
$r=\sqrt{t/T}$. There exist positive
constants $c,C'$ depending on $u,v$ such that $cT\le p_t\le C'T$ for
$m\le t\le n$ and large $T$. Indeed, $p_t\ge N$ and
$p_t-p_{t+1}\le\sqrt{(n-1)/t}$ is uniformly bounded.
The exact identity
\[
 p_t^2-p_{t+1}^2
 =\frac{2\eta_t}{\kappa_0}
  +\left(\frac{\eta_t}{\kappa_0p_{t+1}}\right)^2
\]
therefore gives, uniformly for $m\le t\le n$,
\[
 p_t^2=N^2+2N\sqrt{n-1}\sum_{j=t}^{n-1}j^{-1/2}+O(T).
\]
The squared-increment terms give the $O(T)$ remainder.
Since $N/T\to1-v^2$ and
$T^{-1/2}\sum_{j=t}^{n-1}j^{-1/2}\to 2(v-r)$ uniformly, we have
\[
 \frac{p_t^2}{T^2}\longrightarrow v(1-v^2)(d-4r)
\]
uniformly on this segment, using
$(1-v^2)^2+4v(1-v^2)(v-r)=v(1-v^2)(d-4r)$.
Since $T^{3/2}\kappa_0\to a/[v(1-v^2)]$, it follows from
$M_m=\kappa_0p_m^2$ that $M_m/\sqrt T\to ah$.
Moreover, $\sqrt T\,\ell\to a\sqrt h/(v\sqrt{d-4v})$.
Substitution into the expressions for $\lambda_t$ gives
\begin{equation}\label{eq:cv-lambda-profile}
 \sqrt T\lambda_t\longrightarrow
 \begin{cases}
 \displaystyle\frac{a\sqrt h}{r\sqrt{d-4r}},&u\le r\le v,\\[4pt]
 \displaystyle\frac{a\sqrt h}{v\sqrt{d-4v}},&v\le r\le1.
 \end{cases}
\end{equation}
All denominators are bounded away from zero. Riemann summation gives
\[
 \sum_{t=m}^T\lambda_t^2
 \longrightarrow a^2h\left\{\int_u^v\frac{2\dd r}{r(d-4r)}+\frac1v\right\}
 =a^2hI(u,v).
\]
Substitution in \eqref{eq:cv-potential-bound} proves the claim. Its
right-hand side depends only on scalar coefficients, giving uniformity. 
\end{proof}

\subsection{Trajectories with a common active plane}

The lower construction places the whole trajectory above a single supporting
plane. Two finite-dimensional ingredients provide the required interpolation
and the initialization.

\begin{lemma}[Triangular Gram interpolation]\label{lem:gram}
Let $N\ge2$, let $\lambda_i>0$ for $1\le i\le N$, and define their
tail sums $L_i=\sum_{j=i}^N\lambda_j$.
Assume $\lambda_i\le L_{i+1}$ for $i<N$ and $\lambda_{N-1}=\lambda_N$.
For $0<B_1<L_1$, define
\begin{equation}\label{eq:gram-recursion}
 c_i=\frac{B_i-\lambda_i}{L_{i+1}},\qquad
 B_{i+1}=B_i-\lambda_i c_i\quad(i<N).
\end{equation}
There are unit vectors $\vg_1,\ldots,\vg_N$ in $\R^N$ such that
\begin{equation}\label{eq:gram-common}
 \ip{\vg_i}{\vg_j}=c_i\ (i<j),\qquad
 \vX=\sum_i\lambda_i\vg_i,\quad
 \ip{\vg_i}{\vX}=B_1,\quad \norm{\vX}^2=B_1L_1.
\end{equation}
If the step lengths $\alpha_i>0$ satisfy
$\alpha_i/\lambda_i$ nonincreasing in $i$, then the points $\vz_i$
and scalar values $F_i$, defined for $1\le i\le N$ by
\begin{equation}\label{eq:gram-trajectory}
 \vz_i=\vX-\sum_{j<i}\alpha_j\vg_j,\qquad
 F_i=B_1-\sum_{j<i}\alpha_jc_j
\end{equation}
satisfy $\ip{\vg_j}{\vz_i}\le F_i$ for all $i,j$, with equality when $j\ge i$.
\end{lemma}
\begin{proof}
First check that the recursion permits the construction of the unit
vectors. The recurrence gives, for $i<N$,
\begin{align*}
 B_{i+1}&=\left(1-\frac{\lambda_i}{L_{i+1}}\right)B_i
                 +\frac{\lambda_i^2}{L_{i+1}},\\
 L_{i+1}-B_{i+1}
 &=\frac{L_{i+1}-\lambda_i}{L_{i+1}}(L_i-B_i).
\end{align*}
Consequently $0<B_i\le L_i$ for all $i$, and
$\lambda_{N-1}=L_N=\lambda_N$ gives $B_N=\lambda_N$.

We construct the Gram vectors backwards so that all vectors in a tail
have the same inner product with its weighted sum. Choose a unit vector
$\vg_N$ and set
$\vX_N=\lambda_N\vg_N$. Suppose the tail satisfies
\[
 \vX_{i+1}=\sum_{j>i}\lambda_j\vg_j,\qquad
 \ip{\vg_j}{\vX_{i+1}}=B_{i+1}\quad(j>i),\qquad
 \norm{\vX_{i+1}}^2=B_{i+1}L_{i+1}.
\]
Take a unit vector $\ve_i$ orthogonal to
$\operatorname{span}\{\vg_{i+1},\ldots,\vg_N\}$ (which is possible
in $\R^N$), and define
\[
 \vg_i=\frac{c_i}{B_{i+1}}\vX_{i+1}
       +\sqrt{1-\frac{c_i^2L_{i+1}}{B_{i+1}}}\,\ve_i,
 \qquad \vX_i=\lambda_i\vg_i+\vX_{i+1}.
\]
The square root is real because
\[
 B_{i+1}L_{i+1}-(B_i-\lambda_i)^2=B_i(L_i-B_i)\ge0.
\]
The tail identities give $\norm{\vg_i}=1$ and
$\ip{\vg_i}{\vg_j}=c_i$ for $j>i$. Moreover,
\[
 \ip{\vg_i}{\vX_i}=\lambda_i+c_iL_{i+1}=B_i,\qquad
 \ip{\vg_j}{\vX_i}=\lambda_i c_i+B_{i+1}=B_i\quad(j>i).
\]
Hence $\norm{\vX_i}^2=\sum_{j=i}^N\lambda_j\ip{\vg_j}{\vX_i}=B_iL_i$.
Backward induction, with $\vX=\vX_1$, proves
\eqref{eq:gram-common}, including the case of a zero square root.

\emph{Interpolation along the trajectory.}
To ensure that the plane with slope $\vg_i$ is active at $\vz_i$
in the maximum of all the linear planes, we need
$\ip{\vg_j}{\vz_i}\le F_i$ for every $j$. For this purpose the
recurrence and $L_i=\lambda_i+L_{i+1}$ give
\begin{align*}
 \frac{B_{i+1}}{L_{i+1}}-\frac{B_i}{L_i}
 &=\frac{\lambda_i^2(L_i-B_i)}{L_iL_{i+1}^2}\ge0,\\
 \frac{B_i}{L_i}-c_i
 &=\frac{\lambda_i(L_i-B_i)}{L_iL_{i+1}}\ge0.
\end{align*}
Thus $B_i/L_i$ is nondecreasing and $c_i\le B_i/L_i$.
Fix an old plane $j<i$. For $j\le\ell<i$, telescoping gives
\[
 A_\ell:=\sum_{k=j}^{\ell}\lambda_k
       \bigl(\ip{\vg_j}{\vg_k}-c_k\bigr)
 =B_{\ell+1}-c_jL_{\ell+1}\ge0.
\]
Indeed, the first sum of inner products equals
$B_j-c_jL_{\ell+1}$, since
$B_j=\lambda_j+c_jL_{j+1}$. The sum of the $c_k$ terms is
$B_j-B_{\ell+1}$. The inequality follows from
$c_j\le B_j/L_j\le B_{\ell+1}/L_{\ell+1}$.
Put $r_k=\alpha_k/\lambda_k$. The difference between the value
$F_i$ and the preceding plane is
$F_i-\ip{\vg_j}{\vz_i}
=\sum_{k=j}^{i-1}\alpha_k(\ip{\vg_j}{\vg_k}-c_k)$;
terms with $k<j$ cancel by the Gram identities.
Abel summation now yields
\[
 F_i-\ip{\vg_j}{\vz_i}
 =r_{i-1}A_{i-1}+\sum_{\ell=j}^{i-2}(r_\ell-r_{\ell+1})A_\ell\ge0.
\]
For $j\ge i$, each inner product $\ip{\vg_j}{\vg_k}$ for $k<i$
equals $c_k$, giving equality. 
\end{proof}

\begin{lemma}[Sampling a Gram profile]\label{lem:gram-profile}
Fix $a>0$, $0<u<1$, and a radius $r>0$. Let $\ell$ be a positive,
continuous function on $[u,1]$, with finitely many $C^2$ pieces
and bounded derivatives on each piece, constant near $1$,
with $s\ell(s)$ nondecreasing. Choose an integer $m=u^2T+O(1)$
with $m\ge u^2T$, the first time index of the tail.
Put $N=T-m+1$ and $s_i=\sqrt{(m+i-1)/T}$ for $1\le i\le N$, and set
\[
 b=\frac{\eta_m}{2},\quad \beta=\frac{a}{2u},\qquad
 \lambda_i^0=\frac{\ell(s_i)}{\sqrt T},\qquad
 M(s)=\int_s^1 2z\ell(z)\dd z .
\]
Use the sampled coefficients $\lambda_i^0$ to define
$\lambda_i=(r^2/b)\lambda_i^0/\sum_{j=1}^N\lambda_j^0$, so that
their sum is $L_1=r^2/b$. Take $B_1=b$ and the local step lengths
$\alpha_i=\eta_{m+i-1}$ in Lemma~\ref{lem:gram}.
The construction is valid for sufficiently large $T$, with $\norm{\vX}=r$.
The discrete correlations $c_i$ and values $F_i$ of that lemma satisfy
\begin{equation}
	\begin{aligned}
		\max_{1\le i<N}|Tc_i-c_r(s_i)|&=O(T^{-1}),\\
		\max_{1\le i\le N}|\sqrt T F_i-F_r(s_i)|&=O(T^{-1}),
	\end{aligned}
	\label{eq:correlation-profile}
\end{equation}

\begin{equation}
	\begin{aligned}
		c_r(u)&=\frac{\beta^2}{r^2}-\frac{\ell(u)}{M(u)},
		\qquad&c_r'(s)&=-\frac{\ell'(s)}{M(s)},\\
		F_r(u)&=\beta,&F_r'(s)&=-2ac_r(s).
	\end{aligned}
	\label{eq:gram-limit-ode}
\end{equation}
Here $c_r(s)$ and $F_r(s)$ are the limiting continuous profiles at
radius $r$, determined by \eqref{eq:gram-limit-ode}. Their derivatives
are taken on the smooth pieces, and $c_r$ is constant on the terminal
constant segment of $\ell$. The error constants are uniform in $i$
for the fixed profile and fixed $a,u,r$.
\end{lemma}
\begin{proof}
Write $L_i^0=\sum_{j=i}^N\lambda_j^0$ for the tails of the raw
coefficients, before rescaling. The Riemann sums of
$y\mapsto\ell(\sqrt y)$ on the uniform grid $s_j^2$, of mesh $1/T$,
give $L_i^0/\sqrt T=M(s_i)+O(T^{-1})$ uniformly in $i$.
The constant terminal segment guarantees
$\lambda_i^0\le L_{i+1}^0$ for large $T$ and
$\lambda_{N-1}^0=\lambda_N^0$. Before that segment, the tail sum is
bounded below by a positive multiple of $\sqrt T$. These relations
persist under rescaling. Also, $B_1<L_1$ for large $T$ since $b^2<r^2$,
and $\alpha_i/\lambda_i^0=a/(s_i\ell(s_i))$ is nonincreasing.
The Gram lemma therefore applies and $\norm{\vX}^2=B_1L_1=r^2$.

We obtain the continuous profiles by accumulating the correlation
increments. From \eqref{eq:gram-recursion} and $L_1=r^2/b$,
\[
 c_1=\frac{b^2/r^2-\lambda_1^0/L_1^0}{1-\lambda_1^0/L_1^0}.
\]
Here $s_1=u+O(T^{-1})$, so $Tb^2=\beta^2+O(T^{-1})$ and
$T\lambda_1^0/L_1^0=\ell(u)/M(u)+O(T^{-1})$. Consequently,
\[
 Tc_1=\frac{\beta^2}{r^2}-\frac{\ell(u)}{M(u)}+O(T^{-1}).
\]
Eliminating $B_i$ from the same recurrence gives, for $i<N-1$,
\begin{equation}\label{eq:correlation-difference}
 c_{i+1}-c_i=\delta_i(1-c_i),\qquad
 \delta_i=\frac{\lambda_i^0-\lambda_{i+1}^0}{L_{i+2}^0}.
\end{equation}
To see this, insert $B_i=\lambda_i+c_i L_{i+1}$ into
$B_{i+1}=B_i-\lambda_i c_i$ and use
$L_{i+1}=\lambda_{i+1}+L_{i+2}$. The common rescaling cancels.
The profile $\ell$ is Lipschitz on $[u,1]$.
Whenever its values change between the grid points, $L_{i+2}^0$ is
bounded below by a positive multiple of $\sqrt T$, because $\ell$
is constant near $1$. Since the mesh of $s_i$ is $O(T^{-1})$,
$|\delta_i|=O(T^{-2})$ and $\sum_i|\delta_i|=O(T^{-1})$.
The Gram construction gives $|c_i|\le1$. Summing
\eqref{eq:correlation-difference} first yields
$\max_{i<N}|c_i|=O(T^{-1})$. Using this bound in the same sum gives
\[
 c_i=c_1+\sum_{j<i}\delta_j+O(T^{-2})
 \qquad(1\le i<N).
\]

It remains to identify the sum of the $\delta_j$. On each mesh cell
$[s_j,s_{j+1}]$ before the constant terminal segment,
$L_{j+2}^0/\sqrt T=M(s_j)+O(T^{-1})$, and $M$ stays bounded away
from zero. The integral identity
$\ell(s_{j+1})-\ell(s_j)=\int_{s_j}^{s_{j+1}}\ell'(s)\dd s$
also holds across the piecewise smooth breakpoints. Thus, using
the Lipschitz continuity of $M$,
\[
 T\delta_j
 =-\int_{s_j}^{s_{j+1}}\frac{\ell'(s)}{M(s)}\dd s
   +O(T^{-2}).
\]
The cells contained in the constant terminal segment contribute zero.
Summation and $s_1-u=O(T^{-1})$ give uniformly in $i<N$
\[
 T\sum_{j<i}\delta_j
 =-\int_u^{s_i}\frac{\ell'(s)}{M(s)}\dd s+O(T^{-1}).
\]
Together with the initial value and accumulated increments, this gives
$Tc_i=c_r(s_i)+O(T^{-1})$, where $c_r$ is the profile in
\eqref{eq:gram-limit-ode}.

Finally, \eqref{eq:gram-trajectory} gives
\[
 \sqrt T F_i=\sqrt T b-a\sum_{j<i}\frac{Tc_j}{Ts_j}.
\]
Since $1/(Ts_j)=2(s_{j+1}-s_j)+O(T^{-2})$ and $c_r$ is Lipschitz
on $[u,1]$, substitution and uniform Riemann summation yield
\[
 \sqrt T F_i=\beta-2a\int_u^{s_i}c_r(s)\dd s+O(T^{-1})
 =F_r(s_i)+O(T^{-1}).
\]
This proves both estimates in \eqref{eq:correlation-profile}. 
\end{proof}

\begin{lemma}[Equal-radius initialization]\label{lem:prefix}
Suppose the tail in Lemma~\ref{lem:gram} begins at time $m$, so that
$N=T-m+1$ and $\vx_{m+i-1}=\vz_i$ for $1\le i\le N$.
Use $\alpha_i=\eta_{m+i-1}$, $B_1=b=\eta_m/2$, and suppose $F_i\ge b$.
Put $r^2=bL_1$ and $\rho^2=r^2-b^2$. If $\eta_t b<\rho^2$ for $t<m$,
the tail can be preceded by $m-1$ valid steps such that every prefix point
has norm $r$ and value $b$. All selected gradients have norm at most one,
and all tail points have norm at most $r$. A common unit vector $\vh$ obeys
\begin{equation}\label{eq:common-plane}
 f(\vx_t)=\ip{\vh}{\vx_t}\ge b\quad(1\le t\le T).
\end{equation}
After scaling, this gives a unit-ball instance with hull error at least
$b/\max\{1,r^2\}$.
\end{lemma}
\begin{proof}
Let $\vh=\vg_N$ be the last unit slope of the Gram tail.
By \eqref{eq:gram-common} and \eqref{eq:gram-trajectory},
$\ip{\vh}{\vX}=b$ and $\ip{\vh}{\vz_i}=F_i$, so we aim to make
the plane with slope $\vh$ active also on the prefix. To distinguish
local tail indices from iteration indices, denote the prefix slopes
by $\vq_t$ ($1\le t<m$), retaining $\vg_i$ for the tail slopes.
Choose $\vy_1\perp\vh$ of norm $\rho>0$ in a sufficiently large
Euclidean space, and construct
\begin{equation}\label{eq:equal-radius-prefix}
 \begin{aligned}
 \vx_t&=b\vh+\vy_t,\qquad
 \vq_t=\frac b{\rho^2}\vy_t
       +\sqrt{\frac{2b}{\eta_t}-\frac{b^2}{\rho^2}}\,\ve_t,\\
 \vy_{t+1}&=\vy_t-\eta_t\vq_t,\qquad 1\le t<m .
 \end{aligned}
\end{equation}
Here each $\ve_t$ is a fresh unit vector orthogonal to $\vh$ and to
the span of $\vy_1,\ve_1,\ldots,\ve_{t-1}$; the condition
$\eta_tb<\rho^2$ ensures that the square root is real.
If $\norm{\vy_t}=\rho$, then $\ip{\vq_t}{\vx_t}=b$ and
$\norm{\vq_t}^2=2b/\eta_t\le1$, since $\eta_t\ge\eta_m=2b$.
Expanding $\norm{\vy_t-\eta_t\vq_t}^2$ therefore gives
$\norm{\vy_{t+1}}^2=\rho^2$. Induction gives $\norm{\vx_t}=r$
through $\vx_m$.
The tail can be placed isometrically so that $\vX=b\vh+\vy_m$:
its transverse component has norm $\rho$ by \eqref{eq:gram-common},
and all remaining directions can be taken orthogonal to the
prefix space $\operatorname{span}\{\vy_1,\ve_1,\ldots,\ve_{m-1}\}$.

\emph{Activity of the supporting planes.}
It suffices to check that every linear plane is no higher than the
common plane at each query. Within the prefix, \eqref{eq:equal-radius-prefix}
shows that the projection of $\vy_j$ onto the span occupied at time
$t<j$ is $\vy_t$ multiplied by
$\prod_{s=t}^{j-1}(1-\eta_s b/\rho^2)$. Each factor lies in $(0,1)$.
Hence a future prefix plane has value at most $b$ at a past
prefix point. In the other direction, the inner product of $\vq_t$
with $\vy_{t+1}$ equals $b-\eta_t\norm{\vq_t}^2=-b$, and every
subsequent prefix update multiplies it by a positive factor.
In particular,
\[
 \ip{\vq_t}{\vy_m}
 =-b\prod_{s=t+1}^{m-1}\left(1-\frac{\eta_s b}{\rho^2}\right)\in[-b,0].
\]
This also shows that past prefix planes have value at most $b$
at all later prefix points.

For a tail slope $\vg_i$, set $c_N=1$ so that
$\ip{\vg_i}{\vh}=c_i$ for all $i$. Since $\ip{\vg_i}{\vX}=b$,
its coefficient along $\vy_m$ is $b(1-c_i)/\rho^2\ge0$.
By the placement of the tail, the plane's value at a prefix point is
\[
 bc_i+\frac{b(1-c_i)}{\rho^2}\ip{\vy_m}{\vy_t}\le b.
\]
For the converse comparison, first note that the tail norms do not
increase, since for each $i<N$
\[
 \norm{\vz_{i+1}}^2-\norm{\vz_i}^2
 =\alpha_i^2-2\alpha_iF_i\le0.
\]
Here $\alpha_i\le\eta_m=2b\le 2F_i$.
Decompose the tail point as $\vz_i=F_i\vh+\gamma_i\vy_m+\vzeta_i$,
where $\gamma_i=\ip{\vz_i}{\vy_m}/\rho^2$ and $\vzeta_i$ is
orthogonal to $\vh$ and the prefix space. Since $\norm{\vz_i}\le r$
and $F_i\ge b$,
$\gamma_i^2\rho^2\le r^2-F_i^2\le\rho^2$. Thus $|\gamma_i|\le1$ and
$\ip{\vq_t}{\vz_i}=\gamma_i\ip{\vq_t}{\vy_m}\le b\le F_i$.
Lemma~\ref{lem:gram} handles the comparisons between tail planes.

Define $f$ on the enlarged Euclidean space by
\[
 f(\vx)=\max\bigl\{0,\max_{1\le t<m}\ip{\vq_t}{\vx},
                         \max_{1\le i\le N}\ip{\vg_i}{\vx}\bigr\}.
\]
The prefix maximum is omitted if $m=1$. The function is convex
and globally $1$-Lipschitz, every selected plane is active,
and \eqref{eq:common-plane} holds. The origin minimizes $f$ with
$f(0)=0$; every convex combination of queried points has objective
at least $b$ by the common plane. Projection is inactive for the
constructed updates in the radius-$r$ ball. With $M=\max\{1,r\}$,
divide the points and selected gradients by $M$ and replace $f$ by
$f(M\,\cdot)/M^2$. The scaled points lie in the unit ball, and the
scaled steps are still $\eta_t=a/\sqrt t$. Their hull gap is at least
$b/M^2$, as claimed. 
\end{proof}

\subsection{A matching lower family and scalar minimization}

For $0<u<1/\sqrt3<v<1$, let $d=3v+v^{-1}$ and define the scalar
matching function
\begin{equation}\label{eq:E-definition}
 E(u,v)=\int_u^v\frac{(r-u)(d-6r)}{r^2(d-4r)^2}\dd r.
\end{equation}
For each $u\in(0,1/\sqrt3)$, let $v(u)\in(1/\sqrt3,1)$ solve $E(u,v)=0$
(existence and uniqueness are established in Lemma~\ref{lem:cv-scalar}).
With $v=v(u)$ and $d=3v+v^{-1}$ on this curve, set
\begin{equation}\label{eq:Q-definition}
 Q(u)=\frac{\displaystyle\frac1{v(1-v^2)}+\frac1{d-4u}}{d-2u}.
\end{equation}

\begin{proposition}\label{prop:cv-lower}
For every $0<u<1/\sqrt3$, with $a(u)=2u\sqrt{Q(u)}$,
\begin{equation}\label{eq:cv-parametric-lower}
 \liminf_{T\to\infty}\sqrt T H_T^\cv(a(u))\ge\sqrt{Q(u)}.
\end{equation}
For every fixed $a>0$, also
\begin{equation}\label{eq:cv-orthogonal-lower}
 \liminf_{T\to\infty}\sqrt T H_T^\cv(a)\ge\sqrt{1+a^2/4}.
\end{equation}
\end{proposition}
\begin{proof}
First we construct a tail whose limiting values are all at least
$\sqrt{Q(u)}$. Fix $u$, take $v=v(u)$ and $a=a(u)$, and put
$\beta=\sqrt{Q(u)}=a/(2u)$. In Lemma~\ref{lem:gram-profile},
choose $m=\lceil u^2T\rceil$ and the following profile, with
$s$ denoting the rescaled time variable $\sqrt{t/T}$:
\begin{equation}\label{eq:lower-weights}
 \ell(s)=
 \begin{cases}
 (s\sqrt{d-4s})^{-1},&u\le s<v,\\
 (v\sqrt{d-4v})^{-1},&v\le s\le1 .
 \end{cases}
\end{equation}
The function $s\ell(s)$ is nondecreasing and $\ell$ is constant near
$1$, as required in that lemma. Direct integration, using
$v(d-4v)=1-v^2$, gives
\[
 M(s)=\begin{cases}
 \sqrt{d-4s},&u\le s\le v,\\
 (1-s^2)/(v\sqrt{d-4v}),&v\le s\le1.
 \end{cases}
\]
For radius one, let $c_1(s)$ and $F_1(s)$ denote the continuous
correlation and value profiles in \eqref{eq:gram-limit-ode}.
That identity gives
\[
 c_1(u)=Q(u)-\frac1{u(d-4u)},\qquad
 c_1'(s)=k(s):=\frac{d-6s}{s^2(d-4s)^2}\quad(u<s<v).
\]
On $E(u,v)=0$, Lemma~\ref{lem:cv-scalar} supplies the identity
\begin{equation}\label{eq:Q-integral}
 Q(u)=\frac1{1-v^2}-\int_u^v\frac{2\dd s}{s(d-4s)^2}.
\end{equation}
Since
\[
 k(s)=\frac{2}{s(d-4s)^2}
       -\frac{\mathrm d}{\mathrm ds}\frac1{s(d-4s)},
\]
integration and \eqref{eq:Q-integral} give
\[
 c_1(v)=Q(u)+\int_u^v\frac{2\dd s}{s(d-4s)^2}
                   -\frac1{1-v^2}=0.
\]
Since the correlation profile is constant after $v$, we have
\begin{equation}\label{eq:correlation-return}
 c_1(s)=-\int_{\min\{s,v\}}^v k(z)\dd z,\qquad
 \int_u^1 c_1(s)\dd s=-E(u,v)=0.
\end{equation}
The last equality follows by interchanging the two scalar integrals.
The condition $E(u,v)=0$ forces $u<d/6<v$, since otherwise the
integrand of $E$ would have constant strict sign in the interior.
Thus $k$ changes sign once from positive to negative. The profile
$c_1$ first increases and then decreases to zero at $v$, so it is
positive just before $v$ and zero on the terminal segment.
As its integral is zero, it must start negative, and its negative
values precede its positive values. Consequently,
$\int_u^s c_1(z)\dd z\le0$ for every $s\in[u,1]$, and
$F_1(s)=\beta-2a\int_u^s c_1(z)\dd z\ge\beta$.

To transfer this non-strict continuous bound to the discrete tail,
take any fixed radius $r>1$. The initial correlation in
\eqref{eq:gram-limit-ode} shifts while its derivative stays the same,
so
\[
 \begin{aligned}
 c_r(s)&=c_1(s)-\beta^2(1-r^{-2}),\\
 F_r(s)&=F_1(s)+2a\beta^2(1-r^{-2})(s-u).
 \end{aligned}
\]
This gives strictly negative correlations near $u$ and
$F_r(s)>\beta$ at all $s>u$. By the uniform limits in
\eqref{eq:correlation-profile}, for large $T$ the discrete values
increase from their discrete initial value $b$ near the start of the
tail and exceed $b$
on its remaining part. In detail, fix a short interval on which
$c_r<0$ and use the recurrence for $F_i$ there. On the complement,
the continuous $F_r-\beta$ has a positive minimum.
The prefix condition of Lemma~\ref{lem:prefix} also holds for large
$T$, since $\eta_t b\le ab\to0$ and $r^2-b^2\to r^2>0$.
Applying that lemma and scaling the radius-$r$ construction to
the unit ball gives
\[
 \liminf_{T\to\infty}\sqrt T H_T^\cv(a)\ge\frac{\beta}{r^2}.
\]
Letting $r\downarrow1$ proves \eqref{eq:cv-parametric-lower}.

For the separate lower bound \eqref{eq:cv-orthogonal-lower}, take
$m=\lceil a^2(T+1)/(4+a^2)\rceil$, $N=T-m+1$, and
$b=a/(2\sqrt m)$. In Lemma~\ref{lem:gram} use $\lambda_i=b$
and $B_1=b$, which give $c_i=0$ for all $i<N$.
Thus the slopes may be taken orthonormal, with
$\vX=b\sum_{i=1}^N\vg_i$. All tail values are $b$ and the squared
radius is $r^2=Nb^2\le1$ by the choice of $m$.
Since $\rho^2=(N-1)b^2\to1$ and $\eta_tb\le ab\to0$,
the prefix construction applies for large $T$. The resulting hull error
is at least $b$, with $\sqrt T\,b\to\sqrt{1+a^2/4}$. 
\end{proof}

\begin{lemma}[Scalar matching]\label{lem:cv-scalar}
The root $v(u)$ of $E(u,v)=0$ exists uniquely, is smooth, and is
strictly decreasing on $(0,1/\sqrt3)$.
The function $Q$ has a unique global minimum, attained where
\begin{equation}\label{eq:scalar-optimum}
 u=\frac{v(1-v^2)}2,\qquad v^2=\tau_\ast,\qquad Q(u)=C_\cv^2.
\end{equation}
Moreover $a(u)=2u\sqrt{Q(u)}$ is smooth and strictly increasing, with
endpoint limits $0$ and $\sqrt2$ as $u\downarrow0$ and $u\uparrow1/\sqrt3$,
respectively. Thus it maps $(0,1/\sqrt3)$ bijectively onto $(0,\sqrt2)$.
The identity \eqref{eq:Q-integral} holds on the root curve.
\end{lemma}
\begin{proof}
The scalar differentiations and endpoint calculations are given in
Appendix~\ref{app:scalar}. 
\end{proof}

\begin{proof}[Proof (Theorem~\ref{thm:convex})]
For each $a\in(0,\sqrt2)$, Lemma~\ref{lem:cv-scalar} gives a unique
$u\in(0,1/\sqrt3)$ with $a=a(u)$. Proposition~\ref{prop:cv-lower} then
implies $\Gamma_\cv(a)\ge\sqrt{Q(u)}\ge C_\cv$.
The second inequality is strict unless $u=v(1-v^2)/2$ with
$v^2=\tau_\ast$, at which $2u\sqrt{Q(u)}=a_\cv$ by
\eqref{eq:convex-constants}. For $a\ge\sqrt2$,
\eqref{eq:cv-orthogonal-lower} gives $\Gamma_\cv(a)\ge\sqrt{3/2}>C_\cv$.

For the upper bound, take $v=\sqrt{\tau_\ast}$ and
$u=v(1-v^2)/2$. With $A=1+\tau_\ast+2\tau_\ast^2$,
\eqref{eq:tau-root} gives
\[
 h=A/v,\qquad I(u,v)=\frac1{v(1-v^2)}=\frac1{2u}.
\]
The two terms in \eqref{eq:cv-upper-functional} balance at
$a=(hI)^{-1/2}=a_\cv$, yielding
$U(a_\cv,u,v)=\sqrt{I/h}=C_\cv$.
Proposition~\ref{prop:cv-upper} and convexity prove
\eqref{eq:convex-attainment}. At $a=a_\cv$, the lower limit in
\eqref{eq:cv-parametric-lower} is at least $C_\cv$ and the upper
limit is at most $C_\cv$ by \eqref{eq:convex-attainment}.
This proves convergence and sharpness for the whole hull. 
\end{proof}

The best observed value lies between the hull minimum and the weighted
function value in \eqref{eq:convex-attainment}. It therefore has the same
optimized leading constant in the convex class. Section~\ref{sec:online}
approaches this constant using only the iterates.

\section{The sharp strongly convex hull benchmark}\label{sec:strong}

Set $\vy=\mu(\vx-\vx^\ast)/G$, replace $C$ by $\mu(C-\vx^\ast)/G$,
and replace $f$ by $\mu[f(\vx^\ast+(G/\mu)\vy)-f^\ast]/G^2$.
The transformed objective is $1$-strongly convex, the selected
subgradients are divided by $G$, and the projected subgradient
stepsizes become $a/t$. Thus we work with $\mu=G=1$, $\vx^\ast=0$,
and $f^\ast=0$. For the upper-bound arguments we take $a=1$ and retain
the quadratic contribution of the reported mean, capturing the gain
from strong convexity directly in the output certificate.

\subsection{Projection and a finite energy identity}

We first separate the quadratic part of the objective from its convex
remainder. The effect of projection can then be expressed as an
approximate subgradient inequality for that remainder.

\begin{lemma}[Removing the quadratic part]\label{lem:strong-projection}
Let $h(\vx)=f(\vx)-\norm{\vx}^2/2$, which is convex.
There exist vectors $\vv_t$ and
numbers $e_t\ge0$ such that, for every $\vz\in C$,
\begin{align}
 h(\vz)&\ge h(\vx_t)+\ip{\vv_t}{\vz-\vx_t}-e_t,
       \label{eq:approx-support}\\
 \vv_t&=(t-1)\vx_t-t\vx_{t+1}, \qquad
 \norm{\vv_t}^2\le1-2(t-1)e_t.\label{eq:strong-v}
\end{align}
Moreover,
\begin{equation}\label{eq:strong-distance}
 0\le h(\vx_t)\le\norm{\vx_t}\le1,
 \qquad \norm{\vx_t}\le(t-1)^{-1/2}\quad(t\ge2),
\end{equation}
and
\begin{equation}\label{eq:strong-partial-sum}
 \vH_k:=\sum_{t=1}^k\vv_t=-k\vx_{k+1},
 \qquad \norm{\vH_k}\le\sqrt k.
\end{equation}
\end{lemma}
\begin{proof}
First we control the distances of the queried points.
For $\vx\in C$ and $0<\lambda<1$, constrained optimality at $0$ and
strong convexity give
$0\le f(\lambda\vx)\le\lambda f(\vx)-
\lambda(1-\lambda)\norm{\vx}^2/2$.
Divide by $\lambda$ and let $\lambda\downarrow0$ to obtain
$f(\vx)\ge\norm{\vx}^2/2$. The strong support inequality at $\vx_t$
evaluated at $0$ implies
$\ip{\vg_t}{\vx_t}\ge f(\vx_t)+\norm{\vx_t}^2/2
\ge\norm{\vx_t}^2$, whence $\norm{\vx_t}\le1$ and
\[
 \norm{\vg_t-\vx_t}^2
 =\norm{\vg_t}^2-2\ip{\vg_t}{\vx_t}+\norm{\vx_t}^2\le1.
\]
Since $\vg_t-\vx_t\in\partial h(\vx_t)$, its support at zero gives
$0\le h(\vx_t)\le\ip{\vg_t-\vx_t}{\vx_t}\le\norm{\vx_t}$.
Since projection cannot increase the distance to $0\in C$, we have
$\norm{\vx_{t+1}}^2\le\norm{\vx_t-\vg_t/t}^2$. Expanding and using
$\ip{\vg_t}{\vx_t}\ge\norm{\vx_t}^2$ gives
\[
 \norm{\vx_{t+1}}^2
 \le(1-2/t)\norm{\vx_t}^2+t^{-2}.
\]
At $t=1$, this gives $\norm{\vx_2}^2\le1-\norm{\vx_1}^2\le1$.
For $t\ge2$, the coefficient of $\norm{\vx_t}^2$ is nonnegative.
If $\norm{\vx_t}^2\le1/(t-1)$, the right-hand side is at most
$(1-2/t)/(t-1)+t^{-2}\le1/t$. This proves
\eqref{eq:strong-distance} by induction.

To absorb the projection into the subgradient inequality, use the
normal cone to $C$ at $\vy\in C$,
$N_C(\vy)=\{\vn:\ip{\vn}{\vz-\vy}\le0\text{ for all }\vz\in C\}$.
The optimality condition for the projection in \eqref{eq:algorithm-intro}
with stepsize $1/t$ gives
\[
 \vn_t=t(\vx_t-\vx_{t+1})-\vg_t\in N_C(\vx_{t+1}).
\]
In particular, define $e_t=\ip{\vn_t}{\vx_{t+1}-\vx_t}\ge0$,
where nonnegativity follows by taking $\vz=\vx_t$ in the normal
cone inequality. Set $\vv_t=\vg_t-\vx_t+\vn_t$. For $\vz\in C$,
$\ip{\vn_t}{\vz-\vx_t}\le e_t$, so the support inequality
\[
 \begin{aligned}
 h(\vz)&\ge h(\vx_t)+\ip{\vg_t-\vx_t}{\vz-\vx_t}\\
 &\ge h(\vx_t)+\ip{\vv_t}{\vz-\vx_t}-e_t
 \end{aligned}
\]
implies
\eqref{eq:approx-support}. From
$\vg_t-\vx_t=(t-1)\vx_t-t\vx_{t+1}-\vn_t$, we obtain
\[
 \norm{\vv_t}^2-\norm{\vg_t-\vx_t}^2
 =-2(t-1)e_t-2\ip{\vn_t}{\vx_{t+1}}-\norm{\vn_t}^2
 \le-2(t-1)e_t.
\]
Here $\ip{\vn_t}{\vx_{t+1}}\ge0$ because $0\in C$.
The update gives the expression for $\vv_t$, which telescopes to
\eqref{eq:strong-partial-sum}. 
\end{proof}

\begin{lemma}[Finite energy identity]\label{lem:rank-two}
Let $h,\vv_i,e_i$ be as in Lemma~\ref{lem:strong-projection}.
Choose nonnegative scalar sequences $q_i,z_i$ for $1\le i\le T$,
and define the symmetric kernel $K_{ij}=q_iq_j+z_iz_j$
for $1\le i,j\le T$. Assume
\begin{equation}\label{eq:kernel-discrete}
 K_{i-1,j}\le K_{ij}\quad(2\le i\le j\le T).
\end{equation}
Define the scalar coefficients $U_1=0$, $U_i=(i-1)K_{i-1,i}$
for $2\le i\le T$, and
\[
 A_{ij}=\begin{cases}
 (i-1)(K_{ij}-K_{i-1,j}),&2\le i<j,\\
 0,&\text{otherwise}.
 \end{cases}
\]
Set
\begin{equation}\label{eq:rank-weights}
 w_i=U_i-\sum_{j>i}A_{ij},\qquad D=\sum_iw_i,
 \qquad Q_i=\sum_{j>i}\frac{w_j}{j-1}.
\end{equation}
Sums over iteration indices in this lemma run from $1$ to $T$
subject to the indicated restrictions. If $D>0$, let
\[
 J_w=\frac1D\sum_iw_i h(\vx_i)
       +\frac12\norm{\frac1D\sum_iw_i\vx_i}^2.
\]
Then
\begin{equation}\label{eq:rank-bound}
 DJ_w\le\frac12\sum_iK_{ii}
 +\frac12\left\{
 \norm{\sum_i\frac{Q_i}{\sqrt D}\vv_i}^2
 -\norm{\sum_iq_i\vv_i}^2\right\}.
\end{equation}
The auxiliary coefficients $w_i$ may have either sign.
\end{lemma}
\begin{proof}
The proof combines approximate support inequalities with nonnegative
multipliers to produce the desired signed weights $w_i$.
By summation by parts,
\[
 \lambda_i:=U_i-\sum_{j<i}A_{ji}=\sum_{j<i}K_{ji}\ge0,
 \qquad
 w_i=\lambda_i+\sum_{j<i}A_{ji}-\sum_{j>i}A_{ij}.
\]
Here $\lambda_1=0$; for $i\ge2$, the sum of
$(j-1)(K_{ji}-K_{j-1,i})$ over
$2\le j<i$ equals $(i-1)K_{i-1,i}-\sum_{j<i}K_{ji}$.
Since $h(0)=0$, \eqref{eq:approx-support} yields the inequalities
\[
 \begin{aligned}
 h(\vx_i)&\le\ip{\vv_i}{\vx_i}+e_i,\\
 h(\vx_j)-h(\vx_i)&\le\ip{\vv_j}{\vx_j-\vx_i}+e_j \quad(i<j).
 \end{aligned}
\]
Weight the first by $\lambda_i$ and the second by $A_{ij}$,
which are all nonnegative. Summing gives
\begin{equation}\label{eq:rank-support-sum}
 \sum_iw_i h(\vx_i)
 \le\sum_iU_i\ip{\vv_i}{\vx_i}
       -\sum_{i<j}A_{ij}\ip{\vv_j}{\vx_i}+\sum_iU_i e_i.
\end{equation}
The coefficients of $h(\vx_i)$ in the summed support inequalities
are precisely $w_i$, and those of $e_i$ are $U_i$. The $\vx_1$ terms
vanish because $w_1=U_1=A_{1j}=0$.
We now express the right-hand side as quadratic forms in the
vectors $\vv_i$. Substitute
$\vx_i=-\vH_{i-1}/(i-1)$ from \eqref{eq:strong-partial-sum}
for $i\ge2$.
For $k<j$, the coefficient of $\ip{\vv_j}{\vv_k}$ is
\[
 -\frac{U_j}{j-1}+\sum_{i=k+1}^{j-1}\frac{A_{ij}}{i-1}
 =-K_{j-1,j}+\sum_{i=k+1}^{j-1}(K_{ij}-K_{i-1,j})=-K_{kj}.
\]
Thus the two inner-product sums in \eqref{eq:rank-support-sum} equal
\[
 \frac12\sum_iK_{ii}\norm{\vv_i}^2
 -\frac12\norm{\sum_iq_i\vv_i}^2
 -\frac12\norm{\sum_iz_i\vv_i}^2.
\]
Also, interchanging the sums gives
\[
 \sum_iw_i\vx_i=-\sum_{i=2}^T\frac{w_i}{i-1}\sum_{k<i}\vv_k
                 =-\sum_i Q_i\vv_i.
\]
Adding its squared norm divided by $2D$ gives $DJ_w$ on the
left-hand side. By \eqref{eq:strong-v}, the
remaining coefficient of $e_i$ is $U_i-(i-1)K_{ii}\le0$.
Discard it and the nonpositive square involving $z_i$ to obtain
\eqref{eq:rank-bound}. 
\end{proof}

\subsection{Sharp constants for averaging densities}

\begin{lemma}[Sharp density certificate]\label{lem:density}
Let $q,z:[0,1]\to[0,\infty)$ be twice continuously differentiable, with
\begin{equation}\label{eq:density-normalization}
 q'\le0,\qquad q(1)=0,\qquad \int_0^1q(s)\dd s=\int_0^1z(s)\dd s=1.
\end{equation}
Define the integral profiles $P(s)=\int_0^s q(u)\dd u$ and
$L(s)=\int_s^1 z(u)\dd u$. Suppose
\begin{align}
 q(s)P(s)&=z(s)L(s),
       \label{eq:density-compatibility}\\
 \partial_s\{q(s)q(t)+z(s)z(t)\}&\ge0
       \quad(0\le s\le t\le1).\label{eq:density-monotonicity}
\end{align}
Assume $P(s)<1$ for $s<1$.
Set $\omega(s)=-sq'(s)$, the density determining the averaging weights.
For sufficiently large $T$, the average
\[
 \widehat\vx_T=\frac{\sum_{i=1}^T\omega(i/T)\vx_i}
                    {\sum_{i=1}^T\omega(i/T)}
\]
satisfies, uniformly over instances,
\begin{equation}\label{eq:density-bound}
 f(\widehat\vx_T)\le
 \frac1{2T}\int_0^1(q^2+z^2)\dd s+O(T^{-3/2}).
\end{equation}
The leading constant is sharp:
\begin{equation}\label{eq:sharp-density}
 \lim_{T\to\infty}T\sup_{\mathcal I_T^\scv(1)}f(\widehat\vx_T)
 =\frac12\int_0^1(q^2+z^2)\dd s
\end{equation}
in the normalized strongly convex class.
\end{lemma}
\begin{proof}
\emph{Upper bound via a sampled kernel.}
The assumptions on $q$ imply $\omega\ge0$, and integration by parts
gives $\int_0^1\omega(s)\dd s=1$.
Let $N=T-1$, $s_i=(i-1)/N$, and sample
$q_i=q(s_i)/N$, $z_i=z(s_i)/N$ in Lemma~\ref{lem:rank-two}.
The resulting kernel $K_{ij}$ satisfies \eqref{eq:kernel-discrete}
by \eqref{eq:density-monotonicity}. Denote the auxiliary
weights and tail coefficients by $w_i,D,Q_i$ as in
\eqref{eq:rank-weights}.
Differentiation of \eqref{eq:density-compatibility} gives
\begin{equation}\label{eq:density-derivative}
 z'L=q'P+q^2+z^2.
\end{equation}
For $i\ge2$, the exact weights obey
\[
 \frac{w_i}{i-1}=q_{i-1}q_i+z_{i-1}z_i
 -(q_i-q_{i-1})\sum_{j>i}q_j
 -(z_i-z_{i-1})\sum_{j>i}z_j.
\]
Here $q_{i-1}q_i=q(s_i)^2/N^2+O(N^{-3})$ and
$q_i-q_{i-1}=q'(s_i)/N^2+O(N^{-3})$, with the analogous
relations for $z$. The tail sums of $q_j$ and $z_j$ are
$1-P(s_i)+O(N^{-1})$ and $L(s_i)+O(N^{-1})$, respectively.
These uniform estimates yield
\begin{equation}\label{eq:density-weight-expansion}
 \frac{w_i}{i-1}=\frac1{N^2}\{-q'(s_i)+O(N^{-1})\}.
\end{equation}
Indeed, the leading coefficient is
$q^2+z^2-q'(1-P)-z'L=-q'$ by \eqref{eq:density-derivative}.
Since $w_1=0$, summation gives uniformly in $i$
\begin{equation}\label{eq:density-sampling}
 w_i=\frac{\omega(s_i)}N+O(N^{-2}),\qquad
 D=1+O(N^{-1}),\qquad Q_i=\frac{q(s_i)}N+O(N^{-2}).
\end{equation}
For the last estimate, sum \eqref{eq:density-weight-expansion} over $j>i$
and use $q(1)=0$. Thus $D>0$ for large $T$ and
$\max_i|Q_i/\sqrt D-q_i|=O(N^{-2})$.

Abel summation and \eqref{eq:strong-partial-sum} give
\[
 \norm{\sum_iq_i\vv_i}
 =\norm{\sum_{i<T}(q_i-q_{i+1})\vH_i}
 \le\sqrt T q_1=O(N^{-1/2}).
\]
Meanwhile $\norm{\sum_i(Q_i/\sqrt D-q_i)\vv_i}=O(N^{-1})$,
because $\norm{\vv_i}\le1$. The difference of squares in
\eqref{eq:rank-bound} is therefore $O(N^{-3/2})$, and
\[
 J_w\le\frac1{2D}\sum_iK_{ii}+O(N^{-3/2})
      =\frac1{2N}\int_0^1(q^2+z^2)\dd s+O(N^{-3/2}).
\]

We now replace the auxiliary signed weights by the actual averaging
weights. Let
\[
 \widehat w_i=\frac{\omega(i/T)}{\sum_j\omega(j/T)}.
\]
Equation~\eqref{eq:density-sampling} and smoothness give
$\max_i|w_i/D-\widehat w_i|=O(N^{-2})$.
By \eqref{eq:strong-distance}, $\sum_i\norm{\vx_i}=O(\sqrt T)$.
Although the auxiliary weights may be signed, the estimates
\begin{align*}
 \left|\sum_i(\widehat w_i-w_i/D)h(\vx_i)\right|
 &\le\sum_i|\widehat w_i-w_i/D|\,|h(\vx_i)|,\\
 \left\|\sum_i(\widehat w_i-w_i/D)\vx_i\right\|
 &\le\sum_i|\widehat w_i-w_i/D|\,\norm{\vx_i}
\end{align*}
are both $O(N^{-3/2})$. Both mean vectors have norm $O(N^{-1/2})$
since $\max_i |w_i|/D$ and $\max_i\widehat w_i$ are $O(N^{-1})$.
Their squared norms therefore differ by $O(N^{-2})$.
Jensen's inequality for $h$, applied with the nonnegative weights
$\widehat w_i$, gives
$f(\widehat\vx_T)\le J_w+O(N^{-3/2})$.
Replacing $N=T-1$ by $T$ proves the upper bound
\eqref{eq:density-bound}, with constants depending only on $q,z$.

\emph{Matching lower bound for the given density.}
Fix an integer $N_0\ge1$ and a bounded,
nonnegative, piecewise continuous function $B$ on $[0,1]$.
In dimension $T>N_0+1$, write $\vx=(x_1,\ldots,x_T)$ for a
generic point and $\ve_1,\ldots,\ve_T$ for the standard unit
coordinate vectors. At iteration $t$, let $x_{t,j}$ denote
coordinate $j$ of $\vx_t$. Take $C=[0,\infty)^T$, $\vx_1=0$, and set
$b_j=B(j/T)/T$ for $N_0<j<T$ and $b_j=0$ otherwise. Define
\begin{equation}\label{eq:density-hard-function}
 f_T(\vx)=\frac12\norm{\vx}^2+
 \max\left\{0,\max_{N_0<i\le T}
       \left(\sum_{N_0<j<i}b_jx_j-x_i\right)\right\}.
\end{equation}
This function is $1$-strongly convex and has minimum zero at the
origin. We select zero subgradients in the first $N_0$ queries,
so that $\vx_{N_0+1}=0$. At $t>N_0$, select
$\vg_t=\vx_t+\sum_{N_0<j<t} b_j\ve_j-\ve_t$, the subgradient
corresponding to the linear plane indexed by $t$. Inductively the
coordinates $j\ge t$ of $\vx_t$ vanish, so all planes indexed by
$i\ge t$ have the same nonnegative value. The earlier planes
have smaller or equal values, since $b_j,x_j\ge0$. Thus plane $t$
is active and these oracle choices are valid.

The old coordinates update by
$x_{t+1,j}=[(1-1/t)x_{t,j}-b_j/t]_+$ and the new coordinate
is $x_{t+1,t}=1/t$. Removing the clipping from a coordinate after
its birth gives, for $N_0<j<t$,
\[
 x_{t,j}\ \ge\ \frac1{t-1}
                 -b_j\left(1-\frac j{t-1}\right).
\]
Indeed, the untruncated recursion multiplies $x_{t,j}+b_j$
by $1-1/t$, and the product from time $j+1$ to $t-1$ is $j/(t-1)$.
The old coordinates are nonnegative and cannot exceed
$(1-1/t)x_{t,j}$ after an update, so
\[
 \norm{\vx_{t+1}}^2\le (1-1/t)^2\norm{\vx_t}^2+t^{-2}
 \le \frac1t \quad(t>N_0),
\]
where the final inequality follows by induction from $\vx_{N_0+1}=0$.
For $t>N_0$, the selected subgradient has norm at most
$\norm{\vx_t}+\sqrt{1+\sum_j b_j^2}$.
Consequently a common bound for all queries is
$G_{N_0,T}=1+N_0^{-1/2}+\norm{B}_\infty^2/(2T)$,
where $\norm{B}_\infty=\sup_{s\in[0,1]}|B(s)|$.
Scaling the points and domain by $G_{N_0,T}^{-1}$ and the objective
to $f_T(G_{N_0,T}\,\cdot)/G_{N_0,T}^2$ makes them admissible with
$G=1$ and preserves both strong convexity and the steps $1/t$.

For the average $\widehat\vx_T$ of the original unscaled iterates,
let $\widehat x_{T,j}$ denote its $j$th coordinate.
On every compact interval of $s=j/T$ in $(0,1)$ where $B$ is
continuous, the coordinate lower bound
gives
\[
 T\widehat x_{T,j}\ge
 \int_s^1\omega(v)\{v^{-1}-B(s)(1-s/v)\}\dd v+o(1).
\]
The remainder is uniform on such intervals by Riemann summation,
since $v$ stays bounded away from zero and
$T^{-1}\sum_t\omega(t/T)\to1$. Using $\omega(v)=-vq'(v)$ and
$q(1)=0$, integration by parts gives
\[
 \begin{aligned}
 \int_s^1\frac{\omega(v)}{v}\dd v&=q(s),\\
 \int_s^1\omega(v)(1-s/v)\dd v
 &=-\int_s^1(v-s)q'(v)\dd v=1-P(s).
 \end{aligned}
\]
Thus the limiting coordinate lower bound is
$y(s)=q(s)-(1-P(s))B(s)$.
The $T$th coordinate of every queried point is zero. Plane $T$
in \eqref{eq:density-hard-function} therefore implies
$f_T(\widehat\vx_T)\ge\norm{\widehat\vx_T}^2/2+
\sum_{j<T}b_j\widehat x_{T,j}$.
Substituting the limiting coordinate $y$ gives the expression
\[
 \frac12y^2+By
 =\frac12q^2+qPB-\frac12(1-P^2)B^2.
\]
This quadratic in $B(s)$ is concave for $s<1$ because
$P(s)<1$, and its maximum over $B(s)\ge0$ occurs at
$B(s)=q(s)P(s)/(1-P(s)^2)$. To keep the profile bounded even
near $s=1$, for any $0<\varepsilon<1$, take
$B=qP/(1-P^2)$ on $[0,1-\varepsilon]$ and zero elsewhere.
This choice is bounded, and $y=q/(1+P)\ge0$ on its support.
Because $\widehat x_{T,j}$ and $b_j$ are nonnegative, the coordinate
lower bound can be used in the increasing function
$x\mapsto x^2/2+b_jx$ on $[0,\infty)$.
Restricting the sum to $\delta\le j/T\le1-\varepsilon$ for
$0<\delta<1-\varepsilon$ and taking the lower limit thus gives
\[
 \liminf_{T\to\infty} T\sup_{\mathcal I_T^\scv(1)} f(\widehat\vx_T)
 \ge\frac{1}{(1+N_0^{-1/2})^2}
        \int_\delta^{1-\varepsilon}\frac{q(s)^2}{2(1-P(s)^2)}\dd s .
\]
Compatibility \eqref{eq:density-compatibility} implies
$(P^2+L^2)'=2Pq-2Lz=0$, with $P(0)=0$ and $L(0)=1$.
Since $P<1$ on $[0,1)$, $qP=zL$ and $L^2=1-P^2$ yield
$q^2/(1-P^2)=q^2+z^2$ there. Let $N_0\to\infty$,
$\delta\downarrow0$, and $\varepsilon\downarrow0$.
Together with the uniform upper bound, this proves
\eqref{eq:sharp-density}. 
\end{proof}

The hypotheses for both the sine and quadratic averages can be
verified by the same parametrization. Given
$P=\sin\vartheta$ and $L=\cos\vartheta$, where
$\vartheta:[0,1]\to[0,\pi/2]$ increases from $0$ to $\pi/2$, put
$r(s)=\vartheta'(s)$ and
\[
 q=r\cos\vartheta,\qquad z=r\sin\vartheta.
\]
Then $P'=q$, $L'=-z$, and $qP=zL$.
For $s\le t$, writing $\delta=\vartheta(t)-\vartheta(s)$, the
kernel and its derivative are
\begin{align*}
 q(s)q(t)+z(s)z(t)&=r(s)r(t)\cos\delta,\\
 \partial_s\{q(s)q(t)+z(s)z(t)\}
 &=r(t)\{r'(s)\cos\delta+r(s)^2\sin\delta\}.
\end{align*}
Consequently $r\ge0$ and $r'\ge0$ suffice for
\eqref{eq:density-monotonicity}. For the sine average, take
$\vartheta(s)=\alpha s$, $\alpha=\pi/2$. Then
$q(s)=\alpha\cos(\alpha s)$, $z(s)=\alpha\sin(\alpha s)$, and
$\omega(s)=\alpha^2s\sin(\alpha s)$, giving
\eqref{eq:sine-average}. Here $r=\alpha$ is constant.
For the quadratic average, take
\[
 P(s)=\frac{3s-s^3}{2},\qquad
 L(s)=\frac{(1-s^2)\sqrt{4-s^2}}2,\qquad
 \vartheta(s)=\arcsin P(s).
\]
The identity $P^2+L^2=1$ follows by expanding both sides.
The derivative $r=P'/L=3/\sqrt{4-s^2}$ extends smoothly to $s=1$
and is nondecreasing. Thus $q(s)=\tfrac32(1-s^2)$ and
$z(s)=r(s)P(s)$ are smooth, and $\omega(s)=3s^2$.
Both pairs satisfy $q'\le0$, $q(1)=0$, the normalizations
\eqref{eq:density-normalization}, and $P(s)<1$ for $s<1$.
Lemma~\ref{lem:density} gives the two sharp constants from
$q^2+z^2=r^2$:
\begin{equation}\label{eq:strong-density-values}
 \frac12\int_0^1\alpha^2\dd s=\frac{\pi^2}{8},\qquad
 \frac12\int_0^1\frac{9}{4-s^2}\dd s=\frac{9\log3}{8}.
\end{equation}
Restoring $\mu,G$ proves the sine upper bound
\eqref{eq:sine-bound} and the sharp quadratic limit
\eqref{eq:quadratic-sharp}, including its uniform upper bound.

The value $\pi^2/8$ also has a geometric interpretation. Any compatible
pair traces $P^2+L^2=1$ from $(0,1)$ to $(1,0)$, so the length
$\int_0^1\sqrt{q^2+z^2}\dd s$ is at least $\pi/2$.
Cauchy--Schwarz therefore gives
\[
 \frac12\int_0^1(q^2+z^2)\dd s\ge\frac{\pi^2}{8},
\]
attained by the constant-speed path $\vartheta(s)=\pi s/2$.
We now establish a lower bound over the entire iterate hull for
arbitrary $a$ and prove strictness for $a\ne1$.

\subsection{A separating potential for every convex output}

\begin{proposition}\label{prop:strong-lower}
For every fixed $a>0$, \eqref{eq:strong-universal-lower} holds,
with strict inequality whenever $a\ne1$.
\end{proposition}
\begin{proof}
We construct a family of trajectories that separates the entire
iterate hull from the minimizer. The same separation inequality will
also be used to prove strictness when $a\ne1$.

\emph{Construction and admissibility.}
Fix a radius $0<r<1$, an angle $0<\theta<\pi/2$, and choose an
integer $N\ge1$, the length of the initialization phase, so that
\begin{equation}\label{eq:lower-initialization}
 \frac at\le\min\left\{1,\frac{2r^2}{1+r^2}\right\}\quad(t>N).
\end{equation}
For now these parameters are fixed independently of $T$.
For $T>N+1$, work in $\R^{T+1}$ with coordinates indexed by
$0,1,\ldots,T$ and corresponding unit vectors $\ve_0,\ldots,\ve_T$.
Write $\vx=(x_0,\ldots,x_T)$ for a generic point and $x_{t,j}$
for coordinate $j$ of the $t$th iterate $\vx_t$. Set
\[
 C=[0,\infty)^{T+1},\qquad \beta=\frac{r}{\max\{1,a\}},\qquad
 \vx_1=\beta\ve_0.
\]
To obtain a separating linear inequality for the hull, introduce a
positive vector $\vphi=(\varphi_0,\ldots,\varphi_T)$ with
\[
 \varphi_j=\begin{cases}
 \theta/\beta,&0\le j\le N,\\[2pt]
 \displaystyle\frac{\theta\cos(\theta(N+1)/T)}
                         {\cos(\theta j/T)},&N<j\le T ,
 \end{cases}
\]
The slope of each linear plane will have an entry $-\beta_i$ in
coordinate $i$ and entries $b_j$ in the earlier coordinates.
Choose these coefficients by
\[
 \beta_i=\begin{cases}\beta,&i\le N,\\1,&i>N,\end{cases}
 \qquad
 b_j=\begin{cases}\varphi_{j+1}/\varphi_j-1,&N<j<T,\\
                   0,&\text{otherwise}.
       \end{cases}
\]
Here $1\le i\le T$ and $0\le j\le T$. In particular $b_j\ge0$
and $\varphi_{N+1}=\theta$. Let
\begin{equation}\label{eq:strong-lower-function}
 \vc_i=\sum_{j<i}b_j\ve_j-\beta_i\ve_i,
 \qquad
 f_T(\vx)=\frac12\norm{\vx}^2+
                 \max\{0,\ip{\vc_1}{\vx},\ldots,\ip{\vc_T}{\vx}\}.
\end{equation}
This globally finite function is $1$-strongly convex and has minimum zero at
the origin.

At query $t$, inductively the coordinates of $\vx_t$ with index
$j\ge t$ are zero. Every plane $i\ge t$ has the same nonnegative
value $\sum_{j<t}b_j x_{t,j}$. For $i<t$ its value is
$\sum_{j<i}b_jx_{t,j}-\beta_i x_{t,i}\le\sum_{j<t}b_jx_{t,j}$.
Thus $\vg_t=\vx_t+\vc_t$ is a valid selected subgradient.
Projection onto the orthant acts coordinatewise, giving
\begin{equation}\label{eq:strong-coordinate-update}
 \begin{aligned}
 x_{t+1,j}&=[(1-a/t)x_{t,j}-(a/t)b_j]_+ \quad(j<t),\\
 x_{t+1,t}&=(a/t)\beta_t,\qquad x_{t+1,j}=0\quad(j>t).
 \end{aligned}
\end{equation}
During the first $N$ steps the occupied coordinates have $b_j=0$.
If $a/t\le1$, the new point is a convex combination of $\vx_t$ and
$\beta\ve_t$. If $a/t>1$, only the new coordinate is nonzero and has
size $(a/t)\beta\le r$. In either case the prefix points have norm
at most $r$ and coordinate sum between $\beta$ and $r$.
For $t>N$, \eqref{eq:lower-initialization} gives, with $\eta=a/t$,
\[
 \norm{\vx_{t+1}}^2\le(1-\eta)^2r^2+\eta^2\le r^2.
\]
The second inequality uses
$\eta(1+r^2)\le2r^2$. In addition, writing
$\norm{\vx}_1=\sum_{j=0}^T|x_j|$, we have
$\norm{\vx_{t+1}}_1\le(1-a/t)\norm{\vx_t}_1+a/t\le1$
by induction. Thus all queried points have Euclidean norm at most $r$
and coordinate sum at most one. The selected subgradients have norm at most
\begin{equation}\label{eq:strong-lower-gradient}
 r+\sqrt{1+\sum_jb_j^2}=1+r+o(1),
\end{equation}
because $b_j=O(T^{-1})$ for fixed $\theta$ by Taylor expansion of
the cosine ratio. We will rescale at the end to make the norm bound one.

\emph{A lower bound for the whole hull.}
For $N<j<T$, the choice of $b_j$ gives
$b_j\varphi_j=\varphi_{j+1}-\varphi_j$. Summing and using
$\varphi_{N+1}=\theta$ yields the exact telescoping
identity
\begin{equation}\label{eq:separator-recursion}
 \varphi_t-\sum_{j<t}\varphi_jb_j=\theta\quad(t>N).
\end{equation}
$\ip{\vphi}{\vx_t}\ge\theta$ throughout the prefix
$t\le N+1$, because the coordinate sum is at least $\beta$.
For $t>N$, positivity of $\vphi$, orthant clipping, and
\eqref{eq:separator-recursion} imply
\[
 \ip{\vphi}{\vx_{t+1}}
 \ge(1-a/t)\ip{\vphi}{\vx_t}+(a/t)\theta .
\]
Induction gives the common hull constraint
\begin{equation}\label{eq:hull-separation}
 \ip{\vphi}{\vx}\ge\theta\qquad(\vx\in\mathcal K_T).
\end{equation}

Every point in $\mathcal K_T$ has coordinate $T$ equal to zero. The last
plane in \eqref{eq:strong-lower-function} therefore gives
$f_T(\vx)\ge\norm{\vx}^2/2+\ip{\vb}{\vx}$, where
$\vb=(0,b_1,\ldots,b_{T-1},0)$.
Set $\vy_T=\vphi/T-\vb$. This is the unconstrained minimizer of
the quadratic lower function after subtracting the linear form
$\ip{\vphi}{\vx}/T$. Completing the square in
\[
 T\norm{\vx}^2/2+T\ip{\vb}{\vx}
 =T\norm{\vx-\vy_T}^2/2+\ip{\vphi}{\vx}-T\norm{\vy_T}^2/2
\]
and using \eqref{eq:hull-separation} gives
\begin{equation}\label{eq:separating-square}
 Tf_T(\vx)\ge \theta-\frac1{2T}\norm{\vphi-T\vb}^2+(\ip{\vphi}{\vx}-\theta)+\frac T2\norm{\vx-\vy_T}^2
 \quad(\vx\in\mathcal K_T).
\end{equation}
The last two terms are nonnegative, and their simultaneous
smallness is necessary for the first line to be asymptotically tight.
For fixed $r,\theta,N$, the cosine formula and its Taylor expansion
give, uniformly for the tail indices $N<j<T$,
\[
 \varphi_j=\theta\secant(\theta j/T)+O(T^{-1}),\qquad
 Tb_j=\theta\tan(\theta j/T)+O(T^{-1}).
\]
The fixed prefix and the last coordinate contribute $O(T^{-1})$ to the
normalized squared norm.
Riemann summation gives the limit of the first line of
\eqref{eq:separating-square}:
\begin{equation}\label{eq:strong-lower-scalar}
 \theta-\frac{\theta^2}{2}\int_0^1
          \{\secant(\theta s)-\tan(\theta s)\}^2\dd s .
\end{equation}
Using $\secant^2-\tan^2=1$ and the antiderivative
$2\tan z-z-2\secant z$ of $(\secant z-\tan z)^2$,
the expression equals
\begin{equation}\label{eq:strong-lower-angle}
 A(\theta):=\frac{\theta^2}{2}+\theta(\secant\theta-\tan\theta).
\end{equation}

For fixed $r,\theta$, the bound in \eqref{eq:strong-lower-gradient} is at
most $G_r=1+2r$ for large $T$. Scale the points and feasible set by $G_r^{-1}$
and replace $f_T$ by $f_T(G_r\,\cdot)/G_r^2$. This preserves strong
convexity parameter one and the raw steps $a/t$, and makes the gradient bound
at most one. The resulting admissible instances give
\[
 \liminf_{T\to\infty}T H_T^\scv(a)
 \ge\frac{A(\theta)}{(1+2r)^2}.
\]
Let $r\downarrow0$ and $\theta\uparrow\pi/2$, using
$\secant\theta-\tan\theta=\cos\theta/(1+\sin\theta)\to0$.
This proves \eqref{eq:strong-universal-lower}.

\emph{Strictness for every fixed $a\ne1$.}
We now show directly that the last two terms in
\eqref{eq:separating-square} have a strictly positive combined lower
bound over the entire hull. For any output
$\vx=\sum_{t=1}^T w_t\vx_t$, $\vw\in\mathcal W_T$, define its
two residuals
\begin{equation}\label{eq:strong-remainders}
 S_T=\ip{\vphi}{\vx}-\theta\ge0,\qquad
 D_T=T\norm{\vx-\vy_T}^2.
\end{equation}
Write $\alpha=\pi/2$. The cosine formula for $\varphi_j$ and
the expansion of $b_j$ show that, for fixed $r,\theta,N$,
$T(y_T)_j$ approaches
$\theta\{\secant(\theta s)-\tan(\theta s)\}$ uniformly for
$s=j/T$ in compact subintervals of $(0,1)$. As $\theta\uparrow\alpha$,
this profile converges uniformly on such intervals to
\begin{equation}\label{eq:strong-target-profile}
 Y(s)=\frac{\alpha\cos(\alpha s)}{1+\sin(\alpha s)}>0.
\end{equation}
In the following, all coordinate intervals are fixed inside $(0,1)$
before the limits are taken. For an interval $J$, $|J|$ denotes its length.

First we bound the weight placed on late iterates. Fix $0<q<1$,
choose a closed interval $I\subset(q,1)$ with positive length,
and set $B=\frac12\int_I Y(s)\dd s>0$. For $\theta$ sufficiently
close to $\alpha$, and for all sufficiently large $T$ at fixed
$r,\theta,N$, the preceding profile limit gives
$\sum_{j/T\in I}(y_T)_j\ge B$.
Every iterate with $t\le qT$ has zero coordinates on this block,
while all iterates have coordinate sum at most one. Hence the total
late weight $p=\sum_{t>qT}w_t$ satisfies
\[
 B\le\sum_{j/T\in I}x_j+
       \left|\sum_{j/T\in I}((y_T)_j-x_j)\right|
 \le p+\sqrt{D_T},
\]
where the last inequality uses Cauchy--Schwarz and the fact that
the block contains at most $T$ coordinates. Thus
\begin{equation}\label{eq:positive-late-weight}
 p\ge B-\sqrt{D_T}.
\end{equation}
An output close to the target $\vy_T$ must therefore retain a
positive mass of late iterates.

To study the effect of this mass for $a\ne1$, consider the auxiliary
trajectory $\widetilde\vx_t$, starting at
$\vx_{N+1}$, in which we omit the clipping of the old coordinates
in \eqref{eq:strong-coordinate-update} for $t>N$.
Since $1-a/t\ge0$, induction gives
$\vx_t\ge\widetilde\vx_t$ coordinatewise.
For $N<j<t$ the untruncated coordinate has the exact expression
\[
 \widetilde x_{t,j}
 =\left(\frac aj+b_j\right)
       \prod_{h=j+1}^{t-1}(1-a/h)-b_j .
\]
To see this, add $b_j$ to the scalar recurrence and start
from $x_{j+1,j}=a/j$. For $j/T$ in a fixed compact subinterval
of $(0,1)$, the logarithm of the product is
$-a\log(t/j)+O(1/j)$ uniformly for $j<t\le T$.
This follows from summing
$\log(1-a/h)=-a/h+O(h^{-2})$ and comparing the harmonic sum
with its integral. With $s=j/T$ and $v=t/T$, this product converges
uniformly to $(s/v)^a$. Thus the scaled auxiliary coordinate
approaches
\[
 \left(\frac as+\theta\tan(\theta s)\right)(s/v)^a
 -\theta\tan(\theta s)
\]
for fixed $\theta$, uniformly for $s$ in each fixed compact
subinterval of $(0,1)$ and $s<v\le1$. As $\theta\uparrow\alpha$,
this profile converges uniformly on the same sets to
\begin{equation}\label{eq:strong-coordinate-profile}
 F_a(s,v)=\left(\frac as+\alpha\tan(\alpha s)\right)(s/v)^a
                 -\alpha\tan(\alpha s).
\end{equation}

If $a<1$, choose $q=1/2$ and the corresponding $I,B$ in
\eqref{eq:positive-late-weight}. For $0<s<1/2$ and $v\in[1/2,1]$,
$F_a(s,v)\ge a s^{a-1}-\alpha\tan(\alpha s)$, which tends to
$+\infty$ as $s\downarrow0$. The function $Y$ is bounded near zero.
Choose a closed interval $J\subset(0,1/2)$ of positive length,
sufficiently close to zero that $(B/2)F_a(s,v)>Y(s)+2$ for
$s\in J$, $v\in[1/2,1]$. Equations
\eqref{eq:strong-target-profile} and
\eqref{eq:strong-coordinate-profile}, with the uniform limits just
proved and $\vx_t\ge\widetilde\vx_t$ coordinatewise, imply
\[
 \frac B2\,T x_{t,j}\ge T(y_T)_j+1
 \quad(t>T/2,\ j/T\in J)
\]
for $\theta$ sufficiently close to $\alpha$ and all sufficiently
large $T$ at the fixed parameters. If $D_T\ge B^2/4$, then
$D_T/2\ge B^2/8$. Otherwise \eqref{eq:positive-late-weight}
gives $p>B/2$. By nonnegativity of all iterates,
$Tx_j\ge p\min_{t>T/2}Tx_{t,j}\ge (B/2)\min_{t>T/2}Tx_{t,j}$.
Hence $T(x_j-(y_T)_j)\ge1$ on this block, and
\[
 D_T\ge \frac{\#\{j:j/T\in J\}}{T}\ge\frac{|J|}{2},
\]
where $\#\{\cdot\}$ denotes the number of indices.
Both cases yield the uniform positive gap
\[
 S_T+\frac{D_T}{2}\ge
 \delta_a:=\min\{B^2/8,\ |J|/4\}>0.
\]

If $a>1$, the identity \eqref{eq:separator-recursion} gives,
for the auxiliary trajectory,
\[
 \ip{\vphi}{\widetilde\vx_{t+1}}-\theta
 =(1-a/t)(\ip{\vphi}{\widetilde\vx_t}-\theta)\ge0
 \qquad(t>N).
\]
On the other hand, for $s=1-\varepsilon$,
\[
 \alpha\tan(\alpha s)-\frac{a s^{a-1}}{1-s^a}
 =\frac{a-1}{2}+O(\varepsilon)>0
\]
for sufficiently small $\varepsilon>0$. Here
$\alpha\tan(\alpha(1-\varepsilon))=\varepsilon^{-1}+O(\varepsilon)$,
and expansion of $s^a$ gives the second term as
$\varepsilon^{-1}-(a-1)/2+O(\varepsilon)$.
The strict inequality is equivalent to $F_a(s,1)<0$.
Continuity therefore supplies a closed interval $J\subset(0,1)$
of positive length, $q\in(\max J,1)$, and $\xi>0$ such
that $F_a(s,v)\le-2\xi$ on $J\times[q,1]$.
Choose $I,B$ for this $q$ as in \eqref{eq:positive-late-weight}.
For $\theta$ sufficiently close to $\alpha$ and all sufficiently
large $T$, the uniform coordinate limit gives
$\widetilde x_{t,j}\le-\xi/T$ for every $t>qT$ and
$j/T\in J$, while the actual coordinates are nonnegative.
All other coordinate differences $\vx_t-\widetilde\vx_t$ are
also nonnegative. Since $\varphi_j\ge1$ for these tail coordinates
for large $T$ when $\theta$ is close to $\alpha$, we obtain
\[
 \ip{\vphi}{\vx_t}-\theta
 \ge\ip{\vphi}{\vx_t-\widetilde\vx_t}
 \ge \sigma:=\frac{\xi|J|}{2}\qquad(t>qT).
\]
Using \eqref{eq:hull-separation} also for all earlier iterates,
we have $S_T\ge\sigma p$. Again, either $D_T\ge B^2/4$ or
$p>B/2$ by \eqref{eq:positive-late-weight}. Thus
\[
 S_T+\frac{D_T}{2}\ge
 \delta_a:=\min\{B^2/8,\ \sigma B/2\}>0.
\]

In either case, $\delta_a$ depends only on the fixed $a\ne1$
and the intervals chosen from the limiting profiles. In particular,
it is independent of $r,\theta,N,T$. The estimates hold for
$\theta$ sufficiently close to $\alpha$ and every fixed $r$ and
corresponding $N$, once $T$ is large enough. Combining the residual
bound with \eqref{eq:separating-square} and the scaling by $G_r=1+2r$
gives
\[
 \liminf_{T\to\infty}T H_T^\scv(a)
 \ge\frac{A(\theta)+\delta_a}{(1+2r)^2}.
\]
Since $A(\theta)\to\alpha^2/2$, choose $\theta$ close enough to
$\alpha$ and then $r>0$ small enough that the right-hand side
strictly exceeds $\alpha^2/2$. This proves the strict inequality
in \eqref{eq:strong-universal-lower} for every fixed $a\ne1$. 
\end{proof}

\begin{proof}[Proof (Completion of Theorem~\ref{thm:strong})]
Proposition~\ref{prop:strong-lower} gives the lower bound for each
fixed $a>0$ and strictness for $a\ne1$. The sine average lies in
$\mathcal K_T$ and has the uniform bound \eqref{eq:sine-bound}.
Thus $T H_T^\scv(1)\to\pi^2/8$ and
$\Gamma_\scv(a)>\Gamma_\scv(1)$ for $a\ne1$. This proves uniqueness.
\end{proof}

The strongly convex lower example uses active orthant projection, while the
convex common-plane examples have inactive projection on a ball. Both
constructions apply directly to the instance classes in
\eqref{eq:gamma-convex}--\eqref{eq:gamma-strong} and constrain every output in
the queried-point hull.

\section{Online averages near the hull limits}\label{sec:online}

We now compute the exact worst-case leading constants for the two
full-history online averages. The convex argument combines the
projection potential of Section~\ref{sec:convex} with a common-plane
construction that attains its leading terms.

\subsection{Cubic averaging for convex objectives}

Normalize the convex class as in Section~\ref{sec:convex} and set
$a=1/3$. To construct the suffix potential, fix a cutoff fraction
$0<\theta<1$. For a horizon $n$ let
$m=\lfloor\theta n\rfloor$, $r=m/n$, and take $n$ sufficiently large that
$0<m<n$. For the rescaled time $s\in[r,1]$, define the weight tail
$W_r$, its integral $A_r$, and the multiplier profile $S_r$ by
\begin{equation}\label{eq:cubic-profiles}
 \begin{aligned}
 W_r(s)&=\frac{1-s^4}{1-r^4},\qquad
 A_r(s)=2\int_s^1u^{-1/2}W_r(u)\dd u,\\
 S_r(s)&=\sqrt{A_r(r)}\frac{W_r(s)}{\sqrt{A_r(s)}}.
 \end{aligned}
\end{equation}
The terminal value of $S_r$ is defined by continuity.
These profiles are chosen to make $\kappa_t$ nearly constant in
\eqref{eq:prescribed-bound}.

\begin{lemma}[Cubic suffix potential]\label{lem:cubic-segment}
In Lemma~\ref{lem:cv-potential}, applied to the segment $m+1,\ldots,n$,
use weights $\omega_t=t^3/\sum_{j=m+1}^n j^3$ and
$S_t=S_r(t/n)$ for $m\le t\le n$. Then, with
\begin{equation}\label{eq:cubic-energy}
 B(r)=\int_r^1\frac{W_r(s)^2}{sA_r(s)}\dd s,
\end{equation}
the sequence $S_t$ is nondecreasing and the scalar coefficients satisfy
\begin{align}
 \kappa_t&=\frac{1}{a\sqrt n A_r(r)}+O(n^{-3/2}),\notag\\
 \kappa_{t+1}-\kappa_t&=O(n^{-5/2}),\notag\\
 \sum_{t=m+1}^n\nu_t\eta_tS_t&=\frac{aB(r)}{\sqrt n}
                                      +O(n^{-3/2}).
 \label{eq:cubic-kappa-errors}
\end{align}
The estimates hold uniformly for the indicated indices
($m+1\le t<n$ for the differences).
\end{lemma}
\begin{proof}
First check the monotonicity required by the projection potential,
working initially with $r\le s<1$.
Let $d_r(s)=-W_r'(s)=4s^3/(1-r^4)$. Since
$s^{-1/2}/d_r(s)$ decreases,
\[
 A_r(s)=2\int_s^1\frac{u^{-1/2}}{d_r(u)}d_r(u)W_r(u)\dd u
 \le\frac{s^{-1/2}}{d_r(s)}W_r(s)^2.
\]
Thus
\[
 \frac{S_r'}{S_r}=-\frac{d_r}{W_r}+\frac{s^{-1/2}W_r}{A_r}\ge0.
\]
Differentiating $W_r/S_r=\sqrt{A_r/A_r(r)}$ gives
\begin{equation}\label{eq:cubic-profile-identity}
 S_r(r)=1,\qquad S_r'\ge0,\qquad
 -\left(\frac{W_r}{S_r}\right)'(s)
 =\frac{s^{-1/2}S_r(s)}{A_r(r)}.
\end{equation}

We next estimate the discrete potential coefficients. Let $h=1/n$
be the time mesh, and define
\[
 H_h(s)=(1+h)^2-s^2(s+h)^2,\qquad
 Q_h(s)=\frac{H_h(s)}{H_h(r)S_r(s)}.
\]
For an index $m+1\le t\le n$, put $s=(t-1)h$.
The finite weight tail $W_t$ from Lemma~\ref{lem:cv-potential}
satisfies
\[
 W_t=\frac{\sum_{j=t}^n j^3}{\sum_{j=m+1}^n j^3}
     =\frac{H_h(s)}{H_h(r)}
\]
by the identity $\sum_{j=1}^t j^3=t^2(t+1)^2/4$.
Consequently $Q_h(s)=W_t/S_{t-1}$, and
\eqref{eq:prescribed-coefficients} becomes
\begin{equation}\label{eq:cubic-discrete-profile}
 \nu_t=Q_h(s)-Q_h(s+h),\qquad
 \kappa_t=\frac{\sqrt h}{a}\,
 \frac{\sqrt{s+h}}{S_r(s+h)}
 \frac{Q_h(s)-Q_h(s+h)}h.
\end{equation}
To see smoothness at $s=1$, put $v=\sqrt s$. Integration gives
$A_r(s)=4(8-9v+v^9)/[9(1-r^4)]$.
The factorizations $1-v^8=(1-v)L(v)$ and
$8-9v+v^9=(1-v)^2P(v)$ give
\begin{equation}\label{eq:cubic-endpoint}
 S_r(s)=\frac32\sqrt{\frac{A_r(r)}{1-r^4}}
              \frac{L(v)}{\sqrt{P(v)}},\quad
 L(v)=\sum_{j=0}^7v^j,\quad P(v)=\sum_{j=0}^7(8-j)v^j.
\end{equation}
These functions are smooth and positive, uniformly for $r$ near the fixed
positive number $\theta$ and $r\le s\le1$. The difference quotient in
\eqref{eq:cubic-discrete-profile} is the integral of
$-Q_h'(s+\xi h)$ over $0\le\xi\le1$, so it extends smoothly to $h=0$.
At $h=0$ we have $Q_0=W_r/S_r$. The factor multiplying
$\sqrt h/a$ in \eqref{eq:cubic-discrete-profile} therefore equals
$\sqrt s(-Q_0'(s))/S_r(s)=1/A_r(r)$ by
\eqref{eq:cubic-profile-identity}. Smoothness in $h,s$ implies that
this factor differs from the constant by $O(h)$, and its first
$s$ derivative is $O(h)$, uniformly in the grid cells.
As consecutive values of $s$ differ by $h$, these bounds give
the first two estimates in \eqref{eq:cubic-kappa-errors}.
For the last estimate use $\nu_t\eta_tS_t=\kappa_t(\eta_tS_t)^2$
and Riemann summation, which give
\[
 \sum_{t=m+1}^n(\eta_tS_t)^2
 =a^2\int_r^1\frac{S_r(s)^2}{s}\dd s+O(n^{-1})
 =a^2 A_r(r)B(r)+O(n^{-1}).
\]
This proves the lemma. 
\end{proof}

\begin{proof}[Proof (Convex case of Theorem~\ref{thm:online})]
Let $u$ be the root in Theorem~\ref{thm:online} and set $\theta=u^2$,
$m=\lfloor\theta T\rfloor$. We will show that the full-history
potential is attained asymptotically by joining a spherical prefix
to a Gram tail at this cutoff. Write
\[
 Z_n=\sum_{t=1}^n t^3=\frac{n^2(n+1)^2}{4},\qquad
 p_T=Z_m/Z_T,\qquad \omega_t=t^3/Z_T .
\]
Here $\omega_t$ is the full cubic weight at time $t$, and $p_T$
is the total weight of the prefix $1,\ldots,m$.
The integral defining $A_\theta$ gives
\[
 A_\theta(\theta)=\frac{4(8-9u+u^9)}{9(1-u^8)}.
\]
For $a=1/3$, the equation $16-99u+83u^9=0$ is equivalent to
\begin{equation}\label{eq:cubic-root-match}
 a^2A_\theta(\theta)=2u.
\end{equation}

\emph{Upper bound.}
Apply Lemma~\ref{lem:cv-potential} on the entire trajectory,
with $S_t=1$ for $0\le t\le m$ and $S_t=S_{m/T}(t/T)$ for $t>m$.
On the prefix $\nu_t=\omega_t$ and $\kappa_t=\omega_t/\eta_t$
increases in $t$. On the suffix the coefficients of
Lemma~\ref{lem:cubic-segment} are multiplied by $1-p_T$, which
tends to $1-u^8$. At the junction,
\[
 a\sqrt T\kappa_m\to4u^7,\qquad
 a\sqrt T\kappa_{m+1}\to\frac{1-u^8}{A_\theta(\theta)}.
\]
The latter limit is strictly larger: by \eqref{eq:cubic-root-match},
$A_\theta(\theta)=18u$ and the required inequality reduces to
$73u^8<1$, which holds for $u<1/6$.
As $S_{t-1}=1$ up to and including the junction, the potential terms
in \eqref{eq:prescribed-bound} telescope there to $\kappa_{m+1}$.
Within the suffix, the bounded $S_t$ and the difference estimate
in \eqref{eq:cubic-kappa-errors} give a total positive-increment
contribution of $O(T^{-3/2})$. The prefix's energy contribution satisfies
\[
 \frac{\sqrt T}{2}\sum_{t=1}^m\omega_t\eta_t
 =\frac{a\sqrt T}{2Z_T}\sum_{t=1}^m t^{5/2}
 \longrightarrow\frac{4a}{7}u^7 .
\]
The suffix energy is determined by the last estimate in
\eqref{eq:cubic-kappa-errors}. Convexity and
\eqref{eq:prescribed-bound} thus give the uniform bound
\begin{equation}\label{eq:cubic-exact-bound}
 \limsup_{T\to\infty}\sqrt T\sup_{\mathcal I_T^\cv(1/3)}
 f(\vx_T^{(3)})
 \le\frac{4a}{7}u^7+
 (1-u^8)\left\{\frac1{2aA_\theta(\theta)}+\frac a2 B(\theta)\right\}.
\end{equation}
By the factorization used in \eqref{eq:cubic-endpoint} and the
change of variables $s=v^2$,
\[
 B(\theta)=\frac9{2(1-u^8)}\int_u^1\frac{L(v)^2}{vP(v)}\dd v.
\]
Together with $a^2A_\theta(\theta)=2u$ and $a=1/3$, the
first two terms in \eqref{eq:cubic-exact-bound} sum to
$4u^7/21+(1-u^8)/(12u)=1/(12u)+3u^7/28$.
The remaining integral term is the one in \eqref{eq:C3-definition},
so the upper bound equals $C_3$.

\emph{A Gram tail attaining the energy.}
To match the suffix potential, switch to the time variable
$s=\sqrt{t/T}$ and define, for $u\le s\le1$,
\[
 \ell(s)=\frac{aS_\theta(s^2)}s,\qquad
 M(s)=\int_s^1 2z\ell(z)\dd z,\qquad \beta=\frac a{2u},
 \qquad W(s)=\frac{1-s^8}{1-u^8}.
\]
Here $\ell$ is the profile to be sampled in Lemma~\ref{lem:gram-profile},
$M$ its tail integral, and $W$ the limiting normalized cubic weight
tail on the suffix. Using \eqref{eq:cubic-profile-identity},
$M(s)=a\sqrt{A_\theta(\theta)A_\theta(s^2)}$. Indeed, the
derivative of this expression is $-2aS_\theta(s^2)$ and its value
at $1$ is zero. In particular,
\begin{equation}\label{eq:cubic-tail-match}
 W(s)=\frac{S_\theta(s^2)M(s)}{M(u)},\qquad
 M(u)=aA_\theta(\theta),\qquad \beta M(u)=1.
\end{equation}
$\ell(u)=2\beta$ and $s\ell(s)$ is nondecreasing by
\eqref{eq:cubic-profile-identity}.
The equation $\beta M(u)=1$ identifies the cutoff at which the
potential scale matches the unit-radius Gram construction.

For a fixed truncation point $u<v<1$, let
$\ell_v(s)=\ell(\min\{s,v\})$.
It is constant near $1$, and $s\ell_v(s)$ is nondecreasing, as
required by Lemma~\ref{lem:gram-profile}.
Apply Lemma~\ref{lem:gram-profile} with this truncated profile,
radius one, and the tail beginning at $m+1$.
Write $N=T-m$ for the length of the tail and
$s_i=\sqrt{(m+i)/T}$ for its grid. The construction gives an initial
point $\vX$ of norm one with common inner product $b=\eta_{m+1}/2$,
unit slopes $\vg_i$, and tail points $\vz_i$ as in
\eqref{eq:gram-trajectory}, with local step lengths
$\alpha_i=\eta_{m+i}$. Put $\vh=\vg_N$, so that
$\ip{\vh}{\vz_i}=F_i$ and no tail plane exceeds this value.
Denote the limiting scaled correlations and values by
$c^{(v)}(s),F^{(v)}(s)$ to indicate truncation at $v$;
they solve \eqref{eq:gram-limit-ode} at radius one with $\ell_v$
and $M_v(s)=\int_s^1 2z\ell_v(z)\dd z$ in place of $\ell,M$.
For the untruncated profile define $c,F$ on $[u,1)$ by the corresponding
initial-value equations with $\ell,M$.
By \eqref{eq:cubic-tail-match}, they satisfy
\begin{equation}\label{eq:cubic-tail-ode}
 c(u)=-\beta^2,\quad c'(s)=-\ell'(s)/M(s),\qquad
 F(u)=\beta,\quad F'(s)=-2ac(s).
\end{equation}
Although $c$ can be unbounded near $1$, the limit of $F^{(v)}$ is
uniform on $[u,1]$. To verify this, $\ell_v\to\ell$ uniformly,
$M_v(s)\ge c_0(1-s)$ for some $c_0>0$, and the derivatives
of $\ell_v$ are uniformly bounded on the smooth pieces.
Hence $(c^{(v)})'$ converges to $c'$ away from $1$ and
$|(c^{(v)})'(s)|\le C_0/(1-s)$ for a constant $C_0>0$ independent
of $v$ near $1$. Also $c^{(v)}(u)\to c(u)$.
The integrated expression
\[
 F^{(v)}(s)=\beta-2a\left\{(s-u)c^{(v)}(u)
                     +\int_u^s(s-z)(c^{(v)})'(z)\dd z\right\}
\]
converges uniformly because
$\int_u^1(1-z)|(c^{(v)})'(z)-c'(z)|\dd z\to0$ by dominated
convergence; this also defines $F(1)$ by continuity.
The correlations themselves converge uniformly on compact
subintervals of $[u,1)$.

It remains to check feasibility in the unit ball.
Appendix~\ref{app:cubic-arithmetic} proves the explicit margins
\[
 c(s)\le-3/64\quad(u\le s\le2/7),\qquad
 F(s)-a/(2s)\ge5a/64\quad(2/7\le s\le1).
\]
For fixed $v$ sufficiently close to $1$, the preceding convergence
gives $c^{(v)}<0$ on the first interval and
$F^{(v)}(s)>a/(2s)$ on the second, both with strict uniform margins.
By \eqref{eq:correlation-profile}, for large $T$ the discrete
correlations $c_i$ are negative when $s_i\le2/7$. In this portion
the recurrence $F_{i+1}=F_i-\alpha_i c_i$ makes the values
increase from $F_1=b=\eta_{m+1}/2$. Since $\eta_t$ decreases,
$F_i\ge\eta_{m+i}/2$ there. On the remaining portion the same
inequality follows from the uniform value estimate and
$\sqrt T\,\eta_{m+i}=a/s_i$. Therefore for $i<N$,
\[
 \norm{\vz_{i+1}}^2-\norm{\vz_i}^2
 =\eta_{m+i}^2-2\eta_{m+i}F_i\le0,
\]
and all tail points lie in the unit ball.

\emph{Spherical prefix and common active plane.}
We construct the first $m$ updates so that the common plane with
normal $\vh$ is active and each query is on the unit sphere.
Put $b_t=\eta_t/2$ for $1\le t\le m+1$, choose $\vy_1\perp\vh$
with norm $\sqrt{1-b_1^2}$, and write
$\vx_t=b_t\vh+\vy_t$, where $\vy_t\perp\vh$,
$\norm{\vy_t}^2=1-b_t^2$. Define the prefix slope
\[
 \vq_t=d_t\vh+A_t\vy_t+\xi_t\ve_t,\qquad
 d_t=\frac{b_t-b_{t+1}}{\eta_t},\quad
 A_t=\frac{b_t(1-d_t)}{1-b_t^2}.
\]
Here $\ve_t$ is a fresh unit direction orthogonal to $\vh$ and
$\operatorname{span}\{\vy_1,\ve_1,\ldots,\ve_{t-1}\}$.
The coefficient $\xi_t$ is the nonnegative square root of
$1-d_t^2-A_t^2(1-b_t^2)$, ensuring $\norm{\vq_t}=1$.
This is possible because
$d_t^2+A_t^2(1-b_t^2)\le1/4+1/35<1$, using
$b_t\le1/6$ and $0\le d_t<1/2$. Also $1-\eta_t A_t>0$,
since $\eta_t b_t/(1-b_t^2)\le2/35$.
The update $\vx_{t+1}=\vx_t-\eta_t\vq_t$ satisfies
$\ip{\vh}{\vx_{t+1}}=b_{t+1}$ and
\[
 \ip{\vq_t}{\vx_t}=b_t,\qquad
 \norm{\vx_{t+1}}^2=1-2\eta_tb_t+\eta_t^2=1.
\]
Thus $\vy_{t+1}=(1-\eta_t A_t)\vy_t-\eta_t\xi_t\ve_t$.
Place the Gram tail isometrically at $\vx_{m+1}$,
identifying $\vX=b\vh+\vy_{m+1}$, setting $\vx_{m+i}=\vz_i$,
and taking its other transverse directions orthogonal to $\vh$
and the prefix space $\operatorname{span}\{\vy_1,\ve_1,\ldots,\ve_m\}$.

We check all supporting planes. For $k>t$ in the prefix,
$\ip{\vy_k}{\vy_t}\le\norm{\vy_t}^2$, since successive
projections onto the old space multiply by
$1-\eta_jA_j\in(0,1]$. The future plane consequently has
\[
 \ip{\vq_k}{\vx_t}
 \le d_k b_t+
       b_k(1-d_k)\frac{1-b_t^2}{1-b_k^2}\le b_t,
\]
using $b_k\le b_t$. For an old slope $\vq_t$, let
$\vq_t^\perp=\vq_t-d_t\vh$. Its inner product with
$\vy_{t+1}$ equals
$b_t(1-d_t)-\eta_t(1-d_t^2)
=-b_t(1-d_t)(1+2d_t)\le0$. Its inner products with later
$\vy_k$ remain nonpositive, so $\ip{\vq_t}{\vx_k}\le d_t b_k\le b_k$.

For a tail slope $\vg_i$, let $\gamma_i=\ip{\vg_i}{\vh}$.
Because $\ip{\vg_i}{\vX}=b$, its coefficient along $\vy_{m+1}$
is $b(1-\gamma_i)/(1-b^2)\ge0$. The projection property of
the prefix then gives
\[
 \ip{\vg_i}{\vx_t}
 \le\gamma_i b_t+
       b(1-\gamma_i)\frac{1-b_t^2}{1-b^2}\le b_t
 \qquad(t\le m).
\]
Conversely, the tail point $\vz_i$ has a positive component
along $\vy_{m+1}$ for large $T$. In fact,
\begin{align*}
 \ip{\vX}{\vz_i}
 &=1-b\sum_{j<i}\eta_{m+j}
 \ge1-\frac{a^2(1-u)}u+o(1)>0,\\
 b\ip{\vh}{\vz_i}&=bF_i=O(T^{-1}).
\end{align*}
The first inequality follows by summing the inverse-square-root
steps and using $\sqrt T b\to a/(2u)$. The strict sign follows
from the bracket for $u$. Therefore
$\ip{\vy_{m+1}}{\vz_i}=\ip{\vX}{\vz_i}-bF_i>0$.
The other transverse directions of the tail are orthogonal to
all prefix slopes. Hence $\ip{\vq_t^\perp}{\vz_i}\le0$ and
$\ip{\vq_t}{\vz_i}\le d_t F_i\le F_i$.
Together with the Gram interpolation inequalities, this verifies
the activity of the selected planes and of $\vh$ throughout
the joined trajectory.

Take the unit ball as feasible set and the objective
\[
 f(\vx)=\max\bigl\{0,\max_{1\le t\le m}\ip{\vq_t}{\vx},
                        \max_{1\le i\le N}\ip{\vg_i}{\vx}\bigr\}.
\]
All selected slopes have unit norm, the minimum is zero at the origin,
and the constructed updates are feasible within the ball.
At every queried point $f(\vx_t)=\ip{\vh}{\vx_t}$.
On their convex hull, the support plane gives the reverse of
Jensen's inequality, so $f(\vx_T^{(3)})=\sum_t\omega_t f(\vx_t)$.

\emph{Matching the constant.}
The prefix attains the first term of \eqref{eq:cubic-exact-bound}.
For fixed truncation $v$, use the uniform value limit in
Lemma~\ref{lem:gram-profile} and Riemann summation in the variable
$s_i=\sqrt{(m+i)/T}$ to obtain
\[
 \frac{\sqrt T}{1-p_T}\sum_{i=1}^N\omega_{m+i}F_i
 \longrightarrow\int_u^1 F^{(v)}(s)(-W'(s))\dd s.
\]
Indeed, $-W'(s)=8s^7/(1-u^8)$ is the normalized weight density
on the suffix after the change of variables $t/T=s^2$.
Now let $v\uparrow1$, using uniform convergence to $F$.
Integration by parts, using
\eqref{eq:cubic-tail-match} and \eqref{eq:cubic-tail-ode}, gives
\[
 \int_u^1 F(-W')\dd s
 =\beta-\frac1{M(u)}\int_u^1 2s\ell(s)M(s)c(s)\dd s .
\]
From $c'=-\ell'/M$ and $M'=-2s\ell$,
\[
 (cM^2+\ell M)'=-4s\ell Mc-2s\ell^2 .
\]
$cM^2+\ell M$ equals one at $u$ because $c(u)=-\beta^2$,
$\ell(u)=2\beta$, and $M(u)=1/\beta$ by \eqref{eq:cubic-tail-match},
and tends to zero at $1$, since $M(s)=O(1-s)$ and
$c(s)=O(1+|\log(1-s)|)$. Thus
\[
 \int_u^1 F(-W')\dd s
 =\frac{1+\int_u^1 2s\ell(s)^2\dd s}{2M(u)}
 =\frac1{2aA_\theta(\theta)}+\frac a2 B(\theta).
\]
In the last step we used $\int_u^1 2s\ell(s)^2\dd s
=a^2 A_\theta(\theta)B(\theta)$ by the substitution $s^2$.
For every sufficiently close fixed truncation the construction
exists for all large $T$, giving this lower bound in the limit
$v\uparrow1$ for the worst-case lower limit. The suffix has
limiting weight $1-u^8$. Combining with the prefix matches
\eqref{eq:cubic-exact-bound}. Restoring $R,G$ proves the equality
in \eqref{eq:cubic-main}. 
\end{proof}

\subsection{Quadratic averaging and the online recursions}

\begin{proof}[Proof (Completion of Theorem~\ref{thm:online})]
The sharp density result, applied to $\omega(s)=3s^2$ as in
\eqref{eq:strong-density-values}, gives
\eqref{eq:quadratic-main} and \eqref{eq:quadratic-sharp}.
For the recursions, any weighted running average with weights $i^k$
has the form
$\vx_t^{(k)}=\sum_{i=1}^t i^k\vx_i/\sum_{i=1}^t i^k$.
Separating its last term shows that it satisfies
\[
 \vx_t^{(k)}=\vx_{t-1}^{(k)}
       +\frac{t^k}{\sum_{i=1}^t i^k}
           (\vx_t-\vx_{t-1}^{(k)}).
\]
The identities
$\sum_{i=1}^t i^2=t(t+1)(2t+1)/6$ and
$\sum_{i=1}^t i^3=t^2(t+1)^2/4$ give
\eqref{eq:cubic-recursion}--\eqref{eq:quadratic-recursion}. 
\end{proof}

\section{Numerical experiments}\label{sec:experiments}

At a subgradient-query budget $T$, each reported point $\widehat{\vx}_T$
is formed from the first $T$ iterates of \eqref{eq:algorithm-intro}.
We measure its objective gap $e_T=f(\widehat{\vx}_T)-f^\ast$.

\subsection{Hard instances and sharp hull constants}\label{sec:exp-hard}

\begin{figure}[!htbp]
 \centering
 \includegraphics[width=\linewidth]{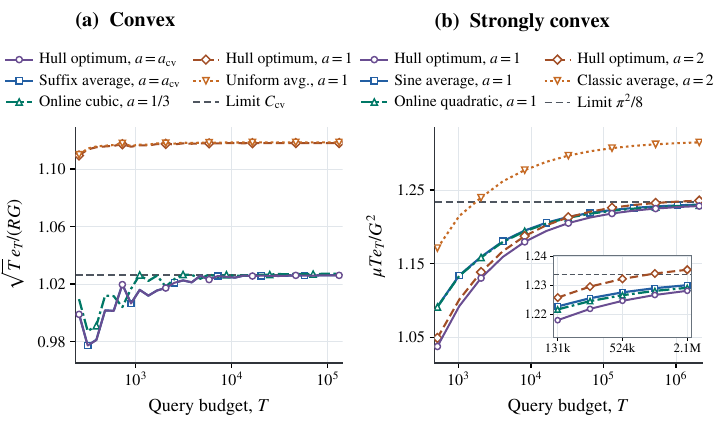}
 \caption{Normalized objective gaps on the convex (a) and strongly convex (b)
 hard instances. The horizontal axes are logarithmic. The dashed reference
 lines mark $C_\cv$ and $\pi^2/8$, respectively. The inset enlarges the two
 hull optima and the proposed averages at large budgets.}
 \label{fig:hard-instances}
\end{figure}

Panels~(a) and~(b) of Figure~\ref{fig:hard-instances} plot
$\sqrt T\,e_T/(RG)$ and $\mu T e_T/G^2$, respectively, against the
budget $T$. The horizontal dashed lines mark $C_\cv$ and $\pi^2/8$,
the sharp asymptotic constants in Theorems~\ref{thm:convex}
and~\ref{thm:strong}. The convex instances use the common-plane
construction of Section~\ref{sec:convex} with $R=G=1$ and a
finite-horizon return condition in the Gram recurrence.  The strongly convex instances use triangular supporting planes as in
Section~\ref{sec:strong} with $\mu=1$.

In general, computing the hull optimum in \eqref{eq:hull-error-intro} requires
the complete trajectory and access to the objective on its hull. The
suffix average at $a=a_\cv$ and the sine average at $a=1$ use
weights prescribed by $t$ and $T$ through
\eqref{eq:convex-optimal-weights} and \eqref{eq:sine-average}. The online
cubic average at $a=1/3$ and quadratic average at $a=1$ use weights
$t^3$ and $t^2$. The recursions \eqref{eq:cubic-recursion} and
\eqref{eq:quadratic-recursion} use the iteration count, the new iterate,
and one running average. The textbook comparisons use the convex uniform
average at $a=1$ and the strongly convex classic average with weights $t-1$
at $a=2$, together with the hull optima at those same coefficients.

At $a=a_\cv$, the hull optimum and suffix average approach $C_\cv$
together. This is the match in Theorem~\ref{thm:convex} between a lower
bound for every convex combination and the upper bound attained by
prescribed weights. The hull optimum at the textbook coefficient $a=1$
stays above $C_\cv$ on its hard family. The separation is already present
after optimization over both hulls, so it shows the benefit of changing
the stepsize coefficient itself. With the nearby coefficient $a=1/3$,
the cubic average tracks the tight curves even when the stopping budget
is chosen during the run, as predicted by Theorem~\ref{thm:online}.

On the strongly convex instances, the sine average and the $a=1$ hull
optimum close in on $\pi^2/8$, illustrating the matching bounds in
Theorem~\ref{thm:strong}. The online quadratic average stays close to
both, so the accuracy gained through the stepsize choice is accessible
with a running average, in agreement with Theorem~\ref{thm:online}.
The inset shows that the $a=2$ hull optimum is higher on the same test
instances and exceeds $\pi^2/8$ at the two largest sampled budgets.
Since each hull optimum gives the best output attainable by averaging
its trajectory, this comparison isolates the improvement due to the
stepsize choice.

At the two optimal coefficients, the small differences at large budgets
between the normalized errors of the prescribed averages and the
instance-specific hull optima illustrate the common limiting value in
\eqref{eq:asymptotic-duality}. The finite-horizon minimax equality
in Theorem~\ref{thm:minimax-duality} concerns worst-case errors over the
full instance class. Likewise, the online curves illustrate the
performance of the running averages on these hull-oriented instances;
their sharp worst-case constants are established in Theorem~\ref{thm:online}.

\FloatBarrier
\subsection{Query efficiency on ball-constrained regression}\label{sec:exp-regression}

We use the same synthetic regression data for ball-constrained versions
of LASSO~\cite{Tibshirani1996} and elastic net~\cite{ZouHastie2005}.

\begin{figure}[!htbp]
 \centering
 \includegraphics[width=\linewidth]{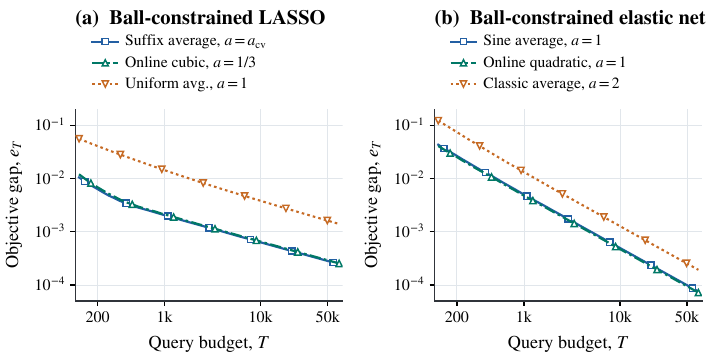}
 \caption{Objective gaps on ball-constrained LASSO (a) and elastic net (b).
 Both axes are logarithmic.}
 \label{fig:sparse-regression}
\end{figure}

Both problems have the form
\[
 \min_{\vx\in\R^{64},\ \norm{\vx}\le r}
 \left\{\norm{\Phi\vx-\vy}^2+\lambda\norm{\vx}_1
                        +\frac{\mu}{2}\norm{\vx}^2\right\},
\]
where $\lambda=0.5$ in both problems, and $\mu=0$ and $0.6$ for LASSO
and elastic net, respectively.
The columns of the Gaussian design $\Phi\in\R^{36\times64}$ are normalized.
The response is $\vy=\Phi\vx_{\rm true}+\bm\varepsilon$, where
$\vx_{\rm true}$ has six nonzero entries and the Gaussian noise has
standard deviation $0.025$. We set $r=0.35\norm{\vx_{\rm true}}\approx0.9008$
and start all runs at $\vx_1=0$. The LASSO runs use $R=2r$ and a common
subgradient bound $G\approx19.9444$. Reference optima are checked against
dual bounds.

Figure~\ref{fig:sparse-regression} shows the unnormalized gap $e_T$
against $T$ on logarithmic axes. The panels use the same stepsize and
weight pairs as the corresponding hard-instance comparisons.

The suffix and online cubic curves stay close to each other and below
the uniform average throughout the sampled LASSO budgets. To reach
$e_T\le 2\times10^{-3}$, the uniform average requires roughly thirty
times as many queries as either proposed average, measured at the first
qualifying budget for each method.

On elastic net, the sine and online quadratic curves also track each
other below the classic average. At the same target gap, the classic
average needs nearly three times as many queries as either proposed
average at the sampled budgets. On both models, the stepsizes and
weights designed for the sharp worst-case bounds thus yield query savings
at a prescribed accuracy. The online recursions retain those savings
with a single running average and a stopping budget chosen during the run.

\FloatBarrier
\section{Conclusions}\label{sec:conclusion}

The iterate hull provides a common worst-case benchmark for all
convex-combination outputs of the projected subgradient method.
Finite-horizon minimax strong duality establishes that, for every finite
horizon and prescribed stepsize sequence, instance-independent weights can
attain the worst-case hull error over the full admissible class. For the inverse-square-root
and harmonic families, we have determined the unique optimal coefficients
$a_\cv$ and $1$ and the corresponding sharp asymptotic constants
$C_\cv$ and $\pi^2/8$, realized by explicit suffix and sine averages.
The cubic and quadratic running averages require neither a prescribed
stopping horizon nor function values, and their exact worst-case leading
constants exceed the optimized limits by only approximately $0.24\%$
and $0.18\%$, respectively. The numerical results illustrate the agreement
with the hull benchmarks and the practical benefit of the resulting
stepsizes and weights on sparse regression: at the sampled budgets, the
proposed rules reach $e_T\le2\times10^{-3}$ with about one thirtieth as
many queries as uniform averaging on LASSO and about one third as many
as the classic average on elastic net. The analysis thereby connects
sharp trajectory-wise benchmarks with implementable averaging rules.

\textbf{Future work:} For a prescribed query budget $T$, Theorem~\ref{thm:minimax-duality}
suggests minimizing the common finite-horizon value $L_T(a)=U_T(a)$
over $a>0$ and characterizing jointly optimal horizon-dependent
stepsize coefficients and instance-independent averaging weights.
Further directions include extending the sharp hull benchmarks to
other prescribed and adaptive stepsize families and to other first-order
methods, such as mirror descent.

\appendix
\section{Scalar matching for the convex benchmark}\label{app:scalar}

We prove Lemma~\ref{lem:cv-scalar} on the parameter region
$0<u<1/\sqrt3<v<1$. Throughout this appendix,
$d=3v+v^{-1}$, $D(r)=d-4r$, and $k(r)=(d-6r)/(r^2D(r)^2)$
for $u\le r\le v$, so that $E(u,v)=\int_u^v(r-u)k(r)\dd r$.
Subscripts $u,v$ on a function of two variables denote partial
derivatives.

\emph{The matching curve.}
We first prove that $E(u,v)=0$ determines a unique smooth root $v(u)$.
With $r$ fixed, $\partial_d k=(8r-d)/(r^2D(r)^3)$. Rearrangement yields
\[
 (r-u)\partial_d k(r)=\frac{\mathrm d}{\mathrm dr}
     \frac{r-u}{rD(r)^2}-\frac{D(u)}{rD(r)^3}.
\]
Since $d'(v)=3-v^{-2}$ and $k(v)=-d'(v)/[vD(v)^2]$, the endpoint
contribution $(v-u)k(v)$ from differentiating $E$ cancels the
term $d'(v)(v-u)/[vD(v)^2]$ from integrating the total derivative.
Thus
\begin{equation}\label{eq:E-v-derivative}
 E_v(u,v)=-(3-v^{-2})(d-4u)
              \int_u^v\frac{\dd r}{r(d-4r)^3}<0.
\end{equation}
To check the endpoint signs, define the auxiliary integral and
its logarithmic and rational terms by
\[
 J=J(u,v)=\int_u^v\frac{2\dd r}{r(d-4r)^2},\qquad
 \Lambda=\log\frac{v(d-4u)}{u(d-4v)},\qquad
 \Delta=\frac1{d-4v}-\frac1{d-4u}.
\]
The partial fraction identity
\[
 \frac{2}{rD(r)^2}=\frac{2}{d^2}\left(\frac1r+\frac4{D(r)}\right)
                       +\frac8{dD(r)^2}
\]
gives
\begin{equation}\label{eq:J-primitive}
 J=\frac2d\Delta+\frac2{d^2}\Lambda,
\end{equation}
Using $k=2/(rD^2)-(1/(rD))'$ and integrating by parts gives
\begin{equation}\label{eq:E-closed}
 E=-\frac{v-u}{v(d-4v)}+\frac\Lambda d+\frac\Delta2-uJ.
\end{equation}
At $v=1/\sqrt3$, $d/6=v$, so the integrand defining $E$ is positive
in the interior.
For fixed $u$, \eqref{eq:J-primitive} and \eqref{eq:E-closed} give
$E(u,v)\sim-(1-u)/[4(1-v)]\to-\infty$ as $v\uparrow1$.
Thus the root $v(u)$ exists uniquely and is smooth.

Its monotonicity follows from the sign change of $k$. On $E=0$, set
$r_0=d/6$. Because the integral of $(r-u)k(r)$ is zero, we have
$u<r_0<v$. The function $k$ is positive before $r_0$ and negative
after $r_0$. Hence
\[
 (r_0-u)\int_u^v k(r)\dd r
 =\int_u^v(r_0-r)k(r)\dd r>0.
\]
It follows that $E_u=-\int_u^v k(r)\dd r<0$. Combining this
with \eqref{eq:E-v-derivative} gives $v'(u)=-E_u/E_v<0$.

For the endpoint limits, fix $0<\varepsilon<1-1/\sqrt3$.
For $1/\sqrt3\le v\le1-\varepsilon$, the expansion
$k(r)=1/(dr^2)+O(r^{-1})$ near zero holds uniformly.
Integrating $(r-u)/(dr^2)$ and the bounded remainder in
$(r-u)k(r)$ gives $E(u,v)=d^{-1}\log(1/u)+O(1)$ as $u\downarrow0$.
The root must therefore satisfy $v(u)\to1$.
At the other endpoint, \eqref{eq:E-v-derivative} remains valid with
$u=1/\sqrt3$ and implies $E(1/\sqrt3,v)<0$ for every $v>1/\sqrt3$.
Continuity and strict decrease in $v$ force $v(u)\to1/\sqrt3$
as $u\uparrow1/\sqrt3$.

\emph{A derivative identity and the step coefficient.}
We next relate the matching curve to the scalar objective $Q(u)$ and
the stepsize coefficient $a(u)$. On $E=0$, eliminate $\Lambda$ using
\eqref{eq:J-primitive} and \eqref{eq:E-closed} to obtain
\[
 J=\frac{\Delta+2(v-u)/\{v(d-4v)\}}{d-2u}.
\]
Since $v(d-4v)=1-v^2$, subtraction from $(1-v^2)^{-1}$ gives
\[
 \frac1{1-v^2}-J
 =\frac{\{v(1-v^2)\}^{-1}+(d-4u)^{-1}}{d-2u}=Q(u).
\]
This proves \eqref{eq:Q-integral} on $E=0$.

Set $\widetilde Q(u,v)=(1-v^2)^{-1}-J(u,v)$, an extension of $Q$
away from the root curve. Differentiating the integral
for $J$ and using $vD(v)=1-v^2$ gives
\[
 \partial_u\widetilde Q=\frac{2}{uD(u)^2},\qquad
 \partial_v\widetilde Q=4d'(v)\int_u^v\frac{\dd r}{rD(r)^3}
                      =-\frac{4E_v}{D(u)}.
\]
Also, integrating $k(r)=2/[rD(r)^2]-\frac{\mathrm d}{\mathrm dr}[1/(rD(r))]$
gives $E_u=-\int_u^v k(r)\dd r=\widetilde Q-1/[uD(u)]$.
Since $E_u+E_vv'(u)=0$ and $Q(u)=\widetilde Q(u,v(u))$, the chain
rule eliminates $v'(u)$:
\begin{equation}
 Q'(u)
 =\frac{2}{uD(u)^2}+\frac{4E_u}{D(u)}=\frac{4Q(u)-2/[uD(u)]}{D(u)}=\frac{2\{2u-v(1-v^2)\}}{u\,v(1-v^2)(d-2u)(d-4u)}.
 \label{eq:Q-derivative}
\end{equation}
In the last equality we used \eqref{eq:Q-definition}.
The expression for $Q'$ in terms of $Q$, together with $a(u)^2=4u^2Q(u)$,
also gives
\[
 \frac{\dd}{\dd u}a(u)^2
 =8uQ(u)+4u^2Q'(u)=\frac{8u}{D(u)}\left\{(d-2u)Q(u)-\frac1{D(u)}\right\}=\frac{8u}{v(1-v^2)D(u)}>0.
\]
Furthermore, $E_u<0$ implies $Q(u)<1/[uD(u)]$.
Thus $0<a(u)^2<4u/D(u)\to0$ as $u\downarrow0$.
At the other endpoint $v(u)\to1/\sqrt3$, so \eqref{eq:Q-definition}
gives $Q(u)\to3/2$ and $a(u)^2\to 4(1/3)(3/2)=2$.
Consequently $a(u)$ is a strictly increasing smooth bijection from
$(0,1/\sqrt3)$ onto $(0,\sqrt2)$.

\emph{The minimum of $Q$.}
It remains to identify the optimal coefficient. Equation
\eqref{eq:Q-derivative} shows that a stationary point has
$u=v(1-v^2)/2$.
With this substitution, write $\tau=v^2$ and $A=1+\tau+2\tau^2$.
Then $D(u)=A/v$, $D(v)=(1-\tau)/v$, and
$\Lambda=\log(2A/(1-\tau)^2)$. Using these in
\eqref{eq:J-primitive}--\eqref{eq:E-closed} yields
\[
 E\left(v(1-v^2)/2,v\right)
 =\frac{v(1+\tau)^2}{(1+3\tau)^2}
       \left\{\log\frac{2A}{(1-\tau)^2}
                   -\frac{1+3\tau}{2(1-\tau)}\right\}.
\]
Hence the stationary point equation is exactly \eqref{eq:tau-root}.
If $g(\tau)$ is the left-hand side of \eqref{eq:tau-root} minus its
right-hand side, then
\[
 g'(\tau)=\frac{1-9\tau^2}{(1+\tau+2\tau^2)(1-\tau)^2}<0
 \quad(\tau>1/3).
\]
Moreover $g(1/3)=\log7-3/2>0$ and
$g(2/3)=\log46-9/2<0$. Conversely, the unique root in $(1/3,2/3)$,
together with $v=\sqrt{\tau_\ast}$ and $u=v(1-v^2)/2$,
lies on $E=0$ in the requisite parameter range and satisfies $Q'(u)=0$.
There is therefore exactly one stationary point.

By the endpoint limits for $v(u)$ and \eqref{eq:Q-definition},
$Q(u)\to\infty$ as $u\downarrow0$ and $Q(u)\to3/2$ as
$u\uparrow1/\sqrt3$. At the stationary point,
\[
 Q(u)=\frac1{2uD(u)}
 =\{(1-\tau_\ast)(1+\tau_\ast+2\tau_\ast^2)\}^{-1}<3/2,
\]
because $(1-\tau)(1+\tau+2\tau^2)$ decreases for $\tau>1/3$
and equals $23/27>2/3$ at $\tau=2/3$.
It is therefore the unique global minimum.
This completes the scalar matching argument.

\section{Scalar feasibility of the cubic construction}\label{app:cubic-arithmetic}

We prove the strict feasibility margins used in the unit-ball
construction for Theorem~\ref{thm:online}. Here $a=1/3$,
and $u$ is the root of $16-99u+83u^9=0$ in $(0,1/6)$.
Indeed, the polynomial is strictly decreasing there and has opposite
signs at the endpoints. Direct substitution of the rational endpoints
of the bracket in Theorem~\ref{thm:online} gives opposite signs as well.
The feasibility proof needs only $0<u<1/6$.
Recall the polynomials $L(s)=\sum_{j=0}^7s^j$ and
$P(s)=\sum_{j=0}^7(8-j)s^j$.
By \eqref{eq:cubic-tail-match} and \eqref{eq:cubic-tail-ode},
the correlation and value profiles satisfy
\[
 \begin{aligned}
 c(u)&=-1/(36u^2),& c'(s)&=-\ell'(s)/M(s),\\
 F(u)&=a/(2u),& F'(s)&=-2ac(s).
 \end{aligned}
\]
We establish
\begin{equation}\label{eq:cubic-feasibility-margins}
 \begin{aligned}
 c(s)&\le-3/64 &&(u\le s\le2/7),\\
 F(s)-a/(2s)&\ge5a/64 &&(2/7\le s\le1).
 \end{aligned}
\end{equation}

With $\theta=u^2$, the factorization in \eqref{eq:cubic-endpoint}
and the formula $M(s)=a\sqrt{A_\theta(\theta)A_\theta(s^2)}$
give, for $K_0=a\sqrt{A_\theta(\theta)/(1-u^8)}>0$,
\[
 \ell(s)=\frac{3K_0}{2}\frac{L(s)}{s\sqrt{P(s)}},\qquad
 M(s)=\frac{2K_0}{3}(1-s)\sqrt{P(s)}.
\]
Differentiating $\ell$ and using $c'=-\ell'/M$ gives
\begin{equation}\label{eq:cubic-c-prime}
 \begin{aligned}
 c'(s)&=-\frac{9N(s)}{8s^2(1-s)P(s)^2},\\
 N(s)&=2sP(s)L'(s)-2P(s)L(s)-sL(s)P'(s).
 \end{aligned}
\end{equation}
The coefficients $n_j$ of $N=\sum_{j=0}^{14}n_js^j$, in ascending order, are
\[
 ( -16,-21,-15,2,30,69,119,180,140,105,75,50,30,15,5 ).
\]
To bound $c'$, expand the following polynomial:
\[
 (1-s)P(s)^2+3N(s)
 =16-15s-12s^2+s^3R(s)+s^{14}(12-s),
\]
where
\[
 R(s)=25+96s+201s^2+340s^3+513s^4+384s^5+287s^6+204s^7+135s^8+80s^9+39s^{10}.
\]
All coefficients of $R$ are positive. For $0<s\le2/3$, the
right-hand side is therefore at least
$16-15s-12s^2\ge2/3>0$. Equation \eqref{eq:cubic-c-prime}
then yields $c'(s)\le3/(8s^2)$ on $[u,2/3]$.
Moreover, the coefficients of $N$ show that $N(s)/s^2$ is strictly
increasing for $s>0$: its negative-power terms have negative
coefficients and its positive-power terms have positive coefficients.
Since $N(2/3)=49352582/4782969>0$, we have $c'(s)<0$ for $2/3\le s<1$.

Integrating the first derivative bound from $u$ to $s\le2/3$
and completing the square gives
\[
 c(s)\le-\frac1{36u^2}+\frac38\left(\frac1u-\frac1s\right)=\frac{81}{64}-\left(\frac1{6u}-\frac98\right)^2-\frac3{8s}\le\frac{81}{64}-\frac3{8s}.
\]
The last expression is at most $-3/64$ on $[u,2/7]$ and at most
$45/64$ on $[u,2/3]$. Since $c$ decreases after $2/3$, we obtain
\[
 c(s)\le-3/64\quad(u\le s\le2/7),\qquad
 c(s)\le45/64\quad(u\le s<1).
\]
Finally, integrating $F'=-2ac$ gives, for $u\le s<1$,
\[
 F(s)\ge a\left(\frac1{2u}-\frac{45}{32}(s-u)\right)
 \ge a\left(\frac1{2u}-\frac{45}{32}(1-u)\right).
\]
The function $v\mapsto 1/(2v)-(45/32)(1-v)$ decreases on
$(0,1/6]$, so $F(s)\ge117a/64$. For $2/7\le s<1$, this gives
\[
 F(s)-\frac a{2s}\ge\frac{117a}{64}-\frac{7a}{4}
 =\frac{5a}{64}.
\]
By \eqref{eq:cubic-c-prime}, $c'(s)=O((1-s)^{-1})$ as
$s\uparrow1$, and hence $c(s)=O(1+|\log(1-s)|)$ is integrable
there. The equation $F'=-2ac$ gives a continuous extension to $1$.
This proves \eqref{eq:cubic-feasibility-margins} also at the endpoint.

\section*{Statements and Declarations}
\ifauthorinformation
\noindent\textbf{Funding.}
This work was supported by National Natural Science Foundation of China
(Grant No. 125B2016, 12631012), and the Academic Excellence Foundation of
BUAA for PhD Students.
\par\medskip
\fi

\ifauthorinformation
\noindent\textbf{Competing Interests.}
The authors declare that there is no conflict of interest.
\par\medskip
\fi

\ifauthorinformation
\noindent\textbf{Data Availability.}
The manuscript has no associated data.
\par\medskip
\fi

\ifauthorinformation
\noindent\textbf{Statement.}
During the preparation of this manuscript, the authors used generative
AI-assisted tools for language refinement, manuscript formatting,
reference organization, and checking and refining proofs. The authors
reviewed and revised all AI-assisted content and take full responsibility
for the accuracy, analyses, and conclusions presented in this work.
\fi



\begin{thebibliography}{99}

\bibitem{Bubeck2015}
Bubeck, S.: Convex optimization: algorithms and complexity.
Found. Trends Mach. Learn. \textbf{8}, 231--357 (2015).
\url{https://doi.org/10.1561/2200000050}

\bibitem{Nesterov2018}
Nesterov, Y.: Lectures on Convex Optimization. Springer, Cham (2018).
\url{https://doi.org/10.1007/978-3-319-91578-4}

\bibitem{ZhuZhangXia2025}
Zhu, Z., Zhang, Y., Xia, Y.: Convergence rate of projected subgradient method
with time-varying step-sizes. Optim. Lett. \textbf{19}, 1027--1031 (2025).
\url{https://doi.org/10.1007/s11590-024-02142-9}

\bibitem{Beck2017}
Beck, A.: First-Order Methods in Optimization. MOS-SIAM Series on
Optimization, vol.~25. Society for Industrial and Applied Mathematics,
Philadelphia (2017).
\url{https://doi.org/10.1137/1.9781611974997}

\bibitem{Nesterov2024}
Nesterov, Y.: Primal subgradient methods with predefined step sizes.
J. Optim. Theory Appl. \textbf{203}, 2083--2115 (2024).
\url{https://doi.org/10.1007/s10957-024-02456-9}

\bibitem{XiaZhangZhu2025}
Xia, Y., Zhang, Y.-H., Zhu, Z.-H.: Lipschitz-free projected subgradient method
with time-varying step-size. J. Oper. Res. Soc. China (2025).
\url{https://doi.org/10.1007/s40305-025-00629-5}

\bibitem{RakhlinShamirSridharan2012}
Rakhlin, A., Shamir, O., Sridharan, K.: Making gradient descent optimal for
strongly convex stochastic optimization. In: Proceedings of the 29th
International Conference on Machine Learning, pp. 1571--1578 (2012).
\url{https://arxiv.org/abs/1109.5647}

\bibitem{ShamirZhang2013}
Shamir, O., Zhang, T.: Stochastic gradient descent for non-smooth optimization:
convergence results and optimal averaging schemes. In: Proceedings of the
30th International Conference on Machine Learning, PMLR \textbf{28}(1),
pp. 71--79 (2013).
\url{https://proceedings.mlr.press/v28/shamir13.html}

\bibitem{LacosteJulienSchmidtBach2012}
Lacoste-Julien, S., Schmidt, M., Bach, F.: A simpler approach to obtaining
an $O(1/t)$ convergence rate for the projected stochastic subgradient method.
arXiv:1212.2002 (2012).
\url{https://doi.org/10.48550/arXiv.1212.2002}

\bibitem{HarveyEtAl2019}
Harvey, N.J.A., Liaw, C., Plan, Y., Randhawa, S.: Tight analyses for non-smooth
stochastic gradient descent. In: Proceedings of the Thirty-Second Conference
on Learning Theory, PMLR \textbf{99}, pp. 1579--1613 (2019).
\url{https://proceedings.mlr.press/v99/harvey19a.html}

\bibitem{GrimmerLi2025}
Grimmer, B., Li, D.: Some primal-dual theory for subgradient methods for
strongly convex optimization. Math. Program. \textbf{214}, 759--788 (2025).
\url{https://doi.org/10.1007/s10107-025-02201-8}

\bibitem{DroriTeboulle2014}
Drori, Y., Teboulle, M.: Performance of first-order methods for smooth convex
minimization: a novel approach. Math. Program. \textbf{145}, 451--482 (2014).
\url{https://doi.org/10.1007/s10107-013-0653-0}

\bibitem{TaylorHendrickxGlineur2017}
Taylor, A.B., Hendrickx, J.M., Glineur, F.: Smooth strongly convex interpolation
and exact worst-case performance of first-order methods. Math. Program.
\textbf{161}, 307--345 (2017).
\url{https://doi.org/10.1007/s10107-016-1009-3}

\bibitem{DroriTaylor2020}
Drori, Y., Taylor, A.B.: Efficient first-order methods for convex minimization:
a constructive approach. Math. Program. \textbf{184}, 183--220 (2020).
\url{https://doi.org/10.1007/s10107-019-01410-2}

\bibitem{TaylorHendrickxGlineurComposite2017}
Taylor, A.B., Hendrickx, J.M., Glineur, F.: Exact worst-case performance of
first-order methods for composite convex optimization. SIAM J. Optim.
\textbf{27}, 1283--1313 (2017).
\url{https://doi.org/10.1137/16M108104X}

\bibitem{LunerGrimmer2026}
Luner, A., Grimmer, B.: On averaging and extrapolation for gradient descent.
Math. Oper. Res. 51(3), 2158--2188 (2026).
\url{https://doi.org/10.1287/moor.2024.0449}

\bibitem{ZamaniGlineur2025}
Zamani, M., Glineur, F.: Exact convergence rate of the last iterate in
subgradient methods. SIAM J. Optim. \textbf{35}, 2182--2201 (2025).
\url{https://doi.org/10.1137/24M1717762}

\bibitem{ZollGrimmer2026}
Zoll, A., Grimmer, B.: A complete characterization of optimal subgradient
methods for Lipschitz convex minimization. arXiv:2607.19240 (2026).
\url{https://doi.org/10.48550/arXiv.2607.19240}

\bibitem{Tibshirani1996}
Tibshirani, R.: Regression shrinkage and selection via the lasso.
J. R. Stat. Soc. B \textbf{58}, 267--288 (1996).
\url{https://doi.org/10.1111/j.2517-6161.1996.tb02080.x}

\bibitem{ZouHastie2005}
Zou, H., Hastie, T.: Regularization and variable selection via the elastic net.
J. R. Stat. Soc. B \textbf{67}, 301--320 (2005).
\url{https://doi.org/10.1111/j.1467-9868.2005.00503.x}


\end{thebibliography}
\end{document}